\documentstyle[12pt,amsfonts]{article}
\newcommand{\mbf}[1]{\mbox{\boldmath$ #1$}}

\newcommand{\be}{\begin{equation}}
\newcommand{\ee}{\end{equation}}
\newcommand{\ba}{\begin{eqnarray}}
\newcommand{\ea}{\end{eqnarray}}

\begin{document}

\begin{center}

{\LARGE\bf  
Electromagnetic pulse propagation in passive 
media by path integral methods}

\vskip 1cm
{\Large Sergei V. Shabanov}\footnote{
electronic mail: {\sf shabanov@phys.ufl.edu\ ; \ \ shabanov@math.ufl.edu}}


\vskip 1cm
{\it Department of Mathematics, University of Florida, Gainesville, FL 32611, USA}

\end{center}

\vskip 1cm

\begin{abstract}

A novel time domain solver of Maxwell's equations in passive
(dispersive and absorbing) media is proposed.
The method is based on the path integral formalism of quantum theory
and entails the use of ({\it i}) the Hamiltonian formalism and 
({\it  ii}) pseudospectral methods (the fast Fourier
transform, in particular) of solving differential equations.
In contrast to finite differencing schemes, the path integral based algorithm 
has no artificial numerical dispersion (dispersive errors),
operates at the Nyquist limit (two grid points per shortest
wavelength in the wavepacket) and exhibits an exponential
convergence as the grid size increases, which, in turn, should 
lead to a higher accuracy. 
The Gauss law holds exactly with no extra computational cost.
Each time step requires $O(N\log_2 N)$
elementary operations where $N$ is the grid size. It can also be
applied to simulations of electromagnetic waves in passive media 
whose properties are time dependent when conventional stationary
(scattering matrix) methods are inapplicable. The stability and
accuracy of the algorithm are investigated in detail.

\end{abstract}

\newpage

\section{Introduction}
\setcounter{equation}0

In this study a time domain solver of Maxwell's equation 
in passive (dispersive and absorbing) media is developed.
The main motivation of this work is to bring methods of 
computational quantum physics into classical
electromagnetic theory. One of the great advantages of 
time domain methods over stationary (scattering
matrix) methods is that a single simulation of the scattering of  a wide band wave packet
can determine basic physical properties of the target 
(e.g., transmission and reflection coefficients)
in the entire frequency band covered by the initial wave packet.
Time domain methods also allow for a unique possibility to observe all immediate effects
on fields caused by the target or by a surrounding passive medium, which greatly facilitates
qualitative understanding of the interaction of an electromagnetic pulse with media and targets.
Another important advantage is that the target geometry (or medium physical properties)
may vary with time and this time dependence cannot be removed by going over to a moving
reference frame. Stationary methods are simply inapplicable to
these kind of problems.

From the computational point of view, the proposed approach is based
on pseudospectral methods. The essential advantages of pseudospectral
algorithms over conventional finite element or finite difference
schemes in solving differential equations
are \cite{psm}: ({\it i}) the exponential versus polynomial rate of convergence
as the grid size (or the basis dimension) increases; ({\it ii})
the absence of dispersive errors and ({\it iii}) efficiency in
numerical calculations.
Time domain algorithms in combination with pseudospectral methods 
have become the state-of-the-art technique in numerical studies of
quantum dynamics by solving the corresponding initial value problem 
for the Schr\"odinger equation (see, e.g., \cite{kosloff}).
A typical algorithm entails an approximate computation
of an object called the path integral (or functional integral)
introduced by Feynman \cite{fey}. Here Maxwell's theory in general
dispersive media is reformulated in the Schr\"odinger (Hamiltonian) formalism. Then
the path integral formalism is  applied to the initial value problem
in  Maxwell's theory in passive media  
to develop a numerical algorithm. 
The main objective of this work is to give a theoretical assessment of the path integral 
based solver of the initial value problem for  
Maxwell's equations. 
Numerical tests and applications will be discussed elsewhere \cite{tests}.

It is shown here that basic principles of the path integral formalism lead
to a true time domain algorithm which indeed enjoys the advantages of pseudospectral
methods. In particular, among the aforementioned features, 
({\it i}) is provided by the use of the fast
Fourier method \cite{fft} as a part of the algorithm, 
when applied to media whose parameters do not have discontinuities in space, ({\it ii}) is a
consequence of Nelson's construction \cite{nelson} of the path integral which
is embedded in our algorithm, ({\it iii}) is due to the fast Fourier
method and some analytical results that speed up numerical computations.
The algorithm has another great advantage over finite difference
schemes: The Gauss law is implemented exactly with no extra
computational cost (Theorem 8.2).
For widely used multi-resonant Lorentz models of passive media,
the algorithm is unitary, meaning that, the energy
of a wave packet is
preserved exactly in dispersive media
with no attenuation (Theorem 6.1). It is also
  unconditionally stable (Theorems 7.1 and 7.3) versus conditionally
stable finite element or finite difference algorithms \cite{psm}
(see also \cite{fd}). 
A possible drawback of the algorithm (to be tested numerically) is that
the use of the fast Fourier method in combination with Nelson's
construction of the path integral might require additional 
computational costs for boundary value problems with complicated
boundary geometry. In our approach, conventional boundary conditions
are not imposed on electromagnetic fields. Targets and medium
interfaces are modeled by discontinuous medium parameters.
The problem arises from 
well known features of the Fourier method \cite{fft}: Aliasing and low convergence
rates for non-smooth functions.
In this study we offer one possible way to 
deal with this problem while keeping the Fourier basis in the
algorithm. Alternative pseudospectral approaches to circumvent
the problem exist and are 
mentioned here, but not discussed in detail.

The basic idea of the path integral 
approach to solving linear, homogeneous, evolutionary 
differential equations (numerically or analytically)
is based on the Hamiltonian formalism. In the framework of 
the Hamiltonian formalism, an original system
of differential equations is transformed to an equivalent 
system of first-order (in time) differential equations
by expanding the original configuration space, that is, by 
going over to a generalized phase space where all time
derivatives, save for the one of highest order, become 
independent variables \cite{arnold}. A generic linear homogeneous 
first-order system can be written in the form
\be
\label{sch}
\partial_t \Psi_t = {\cal H}\Psi_t\ ,\ \ \ \ \ \ \Psi_{t=0} = \Psi_0\ ,
\ee
where $\partial_t$ stands for the partial derivative with respect to time $t$,
a linear operator ${\cal H}$ is called Hamiltonian, while
$\Psi_t$  is called a state vector (or wave function). It is
an element of the generalized
phase space of the system and viewed as a collection (column) of the original variables and
their time derivatives. The generalized phase space is equipped with
an inner product and becomes a Hilbert space. State vectors are
typically vector-valued functions in ${\mathbb R}^3$, and the
Hamiltonian is a differential operator. The choice of the inner
product depends on the problem at hand. One usually requires
that componets of $\Psi_t$ are elements of ${\mathbb L}_2({\mathbb R}^3)$.

In general, upon going over to the Hamiltonian formalism, there might occur
constraints \cite{arnold2,dirac}
\be
\label{constr}
{\cal C}_a \Psi_t = 0\ ,
\ee
with ${\cal C}_a$ being a set of linear operators which do not contain
time derivatives; $a$ enumerates the constraint operators. The constraints must
be preserved in the time evolution which is described by (\ref{sch}). In other
words, the solution is sought in the subspace of the Hilbert
space defined by (\ref{constr}).
Depending on the type of constraints, there are different ways 
of developing the corresponding path integral formalism. In Maxwell's theory,
the constraint is the Gauss law, and it is of the ``first class'' in the Dirac
terminology \cite{dirac}. The characteristic feature of a first class constrained system  is that 
\be
\label{dirac}
[{\cal H}, {\cal C}_a] \sim {\cal C}_a\ ,
\ \ \ \ \ [{\cal C}_a,{\cal C}_b]\sim {\cal
C}_c\ .
\ee
A consequence of (\ref{dirac}) is that if the initial configuration 
$\Psi_0$ satisfies the constraints, then so does the solution of (\ref{sch}).
However, after the projection of the Hilbert space spanned by $\Psi_t$ onto
a finite-dimensional subspace (e.g., a projection on a subspace
associated with a finite spatial grid as is done in Section 3), which is
required for numerical simulations,
the involution condition (\ref{dirac}) can be violated causing
problems 
in simulations. For instance, the Gauss law
is typically violated in any finite differencing approach to
simulations 
of electromagnetic wave packet propagation.
Special efforts have to be made to ensure the transversality of 
the radiation field in Maxwell's theory, which, in turn, complicates
simulation algorithms and increases computational costs 
(e.g., when enforcing the Gauss law 
in finite difference schemes on the grid
via 
the Lagrange multiplier method). It is one of the advantages of the
proposed path integral based algorithm that the Gauss law
can be strictly enforced with no additional computational costs
for generic passive media (Theorem 8.2).

The solution to Eq. (\ref{sch}) is 
\be
\label{evolution}
\Psi_t = \exp\left(t{\cal H}\right) \Psi_0\equiv {\cal U}_t \Psi_0\ ,
\ \ \ t\geq 0\ ,
\ee
assuming that the exponential of ${\cal H}$ exists.
If the Hamiltonian is time dependent then the following replacement 
has to be made in (\ref{evolution})
\be
\nonumber
 \exp\left(t{\cal H}\right) \rightarrow 
T\exp\left(\int_0^t d\tau {\cal H}_\tau\right) = {\cal U}_t\ ,
\ee
where $T\exp$ stands for the time-ordered exponential. The operator
${\cal U}_t$ is defined as the fundamental solution of (\ref{sch}),
$\partial_t {\cal U}_t = {\cal H}_t {\cal U}_t$ with ${\cal U}_{t=0}$
being the identity operator. The fundamental solution has the 
semigroup property,
${\cal U}_{t_1+t_2}= {\cal U}_{t_1}{\cal U}_{t_2}$.
The action of the evolution operator ${\cal U}_t$ on the initial
configuration
can be written via its integral kernel, 
\be
\nonumber
\Psi_t({\bf r}) = 
\int_{{\mathbb R}^3} d{\bf r}^\prime {\cal U}_t({\bf r}, 
{\bf r}^\prime) \Psi_0({\bf r}^\prime) \ . 
\ee
Using the semigroup property of the evolution operator, 
the entire time evolution can be viewed as
consecutive actions of the infinitesimal evolution operator 
${\cal U}_{\Delta t}$, where $\Delta t$ is a time step.
If the kernel of the infinitesimal evolution operator is known, 
then the kernel of the evolution operator can be computed
as the convolution
\be
\label{pi}
{\cal U}_t({\bf r}, {\bf r}^\prime) =\int_{{\mathbb R}^3} 
d{\bf r}_{1}\cdots d{\bf r}_n\,   
{\cal U}_{\Delta t}({\bf r}, {\bf r}_n){\cal U}_{\Delta t}({\bf r}_n,
{\bf r}_{n-1})\cdots {\cal U}_{\Delta t}({\bf r}_1, {\bf r}^\prime)
\ee
with $\Delta t (n+1) = t$. The integration variables can be regarded as points 
${\bf r}_k = {\bf r}( t_k)$, where $t_k = k\Delta t, \ k = 0,1,..., n+1$, 
on a path ${\bf r}(\tau)$ connecting
points ${\bf r}(\tau = t) = {\bf r}$ and ${\bf r}(\tau = 0) = {\bf r}^\prime$.
 In the limit $\Delta t \rightarrow 0$ the convolution (\ref{pi}) can be viewed as a sum over all
paths connecting the initial and final points. 
This is the gist of the Feynman path integral representation
of the fundamental solution of (\ref{sch}). A nontrivial problem is to
find the measure on the space of paths.
For example, if ${\cal H} = \Delta$ (the Laplace operator), it can be
shown that the limit exists, and that the measure coincides with the Wiener
measure which has support in the space of all continuous, but nowhere
differentiable paths (trajectories of the Brownian motion) pinned at the
end points. In quantum mechanics, the problem is more subtle, but can
still be solved \cite{kd}. The existence of the proper measure on the space of
paths opens up an attractive possibility to use Monte-Carlo methods 
of computing the sum over paths which is the gold standard algorithm
in quantum and statistical physics.

However, the present study does not intend to tackle the measure problem 
for the path integral representation of Maxwell's theory, but rather
offers a solution of a more modest problem. Namely, how the
conventional way, outlined above, of deriving the path integral from the original
differential equation can be used to obtain an algorithm for numerical
simulations of the convolution (\ref{pi}) for a small, but finite
$\Delta t$. Similar ideas have been
applied to non-dispersive and/or random media as well as to scattering
problems and waveguides \cite{texas}. Our approach applies to general
passive media and goes beyond the eikonal approximation of
geometric optics and/or the diffraction theory used in earlier works
on path integrals in electromagnetic theory. The results obtained here
are believed to be useful for further development of path integral
methods in theoretical and numerical studies of propagation of electromagnetic 
wave packets in passive media.

The idea of numerical simulations follows from (\ref{pi}) 
rather straightforwardly, namely,
\be
\label{te}
\Psi_{t +\Delta t} = {\cal U}_{\Delta t} \Psi_t\ .
\ee
Thus, finding a state of the system in a sequential moment of time
amounts to computing the action of the exponential of a differential
operator
${\cal H}$ on the state at the preceding moment of time.
Theoretically, it is sufficient to know ${\cal U}_{\Delta t}$ 
up to $(\Delta t)^2$. The limit $\Delta t \rightarrow 0$ in (\ref{pi}) 
would not change if we replace the exact infinitesimal
evolution operator kernel by such an approximation. In numerical simulations,
the limit is never achieved. 
Therefore a higher 
accuracy is required to make errors small.  
Note that the errors are accumulated as more iterations (\ref{te})
are taken. An expansion of $\exp(\Delta t {\cal H})$ into the power series up to some desired
order is known to produce unstable algorithms.
Yet another obvious drawback  is the lack 
of unitarity of the time evolution, that is, if the Hamiltonian is 
skew-symmetric (anti-Hermitian, if a complex phase space is used),
${\cal H}^* = -{\cal H}$, then ${\cal U}_t^* {\cal U}_t = 1$.
In the Maxwell theory, as
we shall see, the squared norm of $\Psi_t$ with respect to the
 ${\mathbb L}_2({\mathbb R}^3)$ scalar product $(\Psi_1, \Psi_2) = \int d{\bf r}
\Psi_1^*\Psi_2$ 
is proportional to the electromagnetic
energy of the system. Consequently, for non-absorbing media
the unitarity of the time evolution is required in simulations to
provide the energy conservation. 

We shall apply Nelson's method of obtaining the path integral 
representation of the fundamental solution of Maxwell's equations
for passive media. It is
based on the Kato-Trotter product formula for the exponential of a sum of two
noncommuting operators and the use of the Fourier basis to compute
exponentials of differential operators. 
Actually, in practical applications, a simpler
version, known as the Lie-Trotter product formula, is used (see the textbooks
\cite{r} for details and references therein). 
In computational quantum mechanics this is also known as the split
operator method. It allows one to keep the differential operators in
the exponential, and thereby, ensures the correct dispersion relation
of simulated electromagnetic waves. It will be shown that there exists
a particular realization of this idea in which the Gauss law holds
exactly in simulations. In general,
the Gauss law can be enforced
by the projection operator formalism developed for the path integral 
representation of constrained dynamical systems (for a review see \cite{pr}
and references therein,  a numerical application to constrained wave
packet propagation can be found in \cite{f}).
The idea is to replace the Hamiltonian by its projection 
on the subspace (\ref{dirac}). If ${\cal P}$ is the projection
operator, that is, ${\cal C}_a{\cal P}\Psi = 0$ for any $\Psi$,
${\cal P}^2 = {\cal P}$ and ${\cal P}^* = {\cal P}$, then 
${\cal H}$ is replaced by ${\cal P}{\cal H}{\cal P}$.
In  Maxwell's theory the projection 
can be implemented in our algorithm with no extra computational costs.
A significant difference from the quantum mechanical
case is that the Hamiltonian ${\cal H}$ (or its projection) is not normal, that is,
it does not commute with its adjoint. This feature complicates
the stability analysis because the von Neumann criteria is no longer
sufficient for stability, while still being necessary \cite{richt}.
Nevertheless, the stability, accuracy and convergence analysis of the
algorithm can be carried out in rather general settings.

Since time domain simulations are performed on finite lattices,
there is always a moment of time when the simulated signal first
reaches the lattice boundary. One typically uses lattices with
periodic boundary conditions. So, the pulse would appear on the
other side of the lattice interfering with itself, thus leading
to totally disastrous results for simulations. The problem
is usually solved by introducing absorbing boundary conditions (see,
e.g., \cite{cond} (for quantum mechanics) and
\cite{abc} (for electrodynamics)).
It is convenient
to set a conducting layer at the grid boundary whose conductivity
is chosen so that it neither transmits nor reflects within the 
designated accuracy in the frequency domain of the initial pulse. 
In Appendix
we briefly describe how such a conducting layer can be 
obtained.

\section{Maxwell theory in the Hamiltonian formalism}
\setcounter{equation}0

Dynamics of electromagnetic waves in continuous media is governed
by Maxwell's equations
\ba
\label{m1}
\partial_t{\bf D}_t &=& c
\mbf{\nabla}\times {\bf H}_t\\
\label{m2}
\partial_t{\bf B}_t &=& -c
\mbf{\nabla}\times {\bf E}_t\ ,
\ea
where $c$ is the speed of light, 
boldface letters denote three-vector fields in ${\mathbb R}^3$ whose
spatial arguments are suppressed and the time dependence is indicated by a subscript.
No external currents and
charges (antennas) are included in this study. 
However, the formalism being
developed is readily generalized to the case when external time
dependent sources are present. 
The electric and magnetic induction vectors, ${\bf D}_t$ and ${\bf B}_t$, respectively,
 are subject to the constraints (the Gauss law)
\be
\label{gauss}
\mbf{\nabla}\cdot{\bf D}_t = 
\mbf{\nabla}\cdot {\bf B}_t = 0\ .
\ee
In linear response theory, assumed through out the paper,
the electric induction is related to the electric field as \cite{landau}
\be
\label{de}
{\bf D}_t = {\bf E}_t + \int^t_{-\infty} d\tau\, \chi_{t-\tau}^e\, 
{\bf E}_{\tau} \equiv {\bf E}_t + {\bf P}_t\ ,
\ee
where $\chi_t^{e}$ is an electric response function of the medium
and ${\bf P}_t$ is the medium polarization vector. A similar relation
can be written for the magnetic field and induction, ${\bf B}_t = {\bf H}_t + {\bf M}_t$,
where magnetization ${\bf M}_t$ is determined by the applied magnetic field and
the magnetic response function of the medium. 

The relation between inductions and fields must be causal, meaning that
the response of the medium, ${\bf P}_t$ and ${\bf M}_t$, can only
depend on fields applied to the medium prior to the current time $t$, 
(e.g., $\chi_t^e = 0$ for $t<0$) \cite{landau}. A natural way to ensure the causality
is to require that the response function satisfies a differential
equation. In other words, the response function is assumed to be the
fundamental solution of some time evolution differential equation.
This differential equation can be obtained from a particular physical
model of the medium in question. A popular model is the multi-resonant Lorentz model.
Let $\tilde{\bf D}_\omega$ and
$\tilde{\bf E}_\omega$ be the Fourier transforms of the electric
induction and field. 
Then from (\ref{de}) it follows that
$\tilde{\bf D}_\omega = \varepsilon_\omega \tilde{\bf E} _\omega$. The dielectric constant in
the Lorentz model has the form
\be
\label{lm}
\varepsilon_\omega = 1 +\sum_{a=1}^N \frac{\omega_{pa}^2}{\omega_a^2 -\omega^2 -2i\gamma_a\omega}\ ,
\ee
and ${\bf M}_t={\bf 0}$.
The physical meaning of the model is transparent. The medium is assumed to be made of 
$N$ sorts of  damped harmonic oscillators
with resonant frequencies $\omega_a$ and damping coefficients
$\gamma_a$. 
Parameters $\omega_{pa}$, called the plasma frequencies, 
are proportional to coupling constants of the
oscillators to the external electric field (the electric dipole
coupling) and also depend
on the density of oscillators of the sort $a$. 
The density may vary in space.
So $\omega_{pa}$ are assumed to be functions of 
spatial coordinates. In an empty space, $\omega_{pa} = 0$.
If the resonant frequency is zero,
the one-resonant Lorentz model is equivalent to the Drude model of
metals \cite{landau}. 
In the Lorentz medium the magnetic response function is zero, while the electric response function 
can easily be found by taking the Fourier transform of (\ref{lm}). Its explicit form is omitted here
because it will not be used. 
The medium polarization is determined by a set of second-order differential 
equations
\be
\label{p}
\partial_t^2 {\bf P}_t^a + 2\gamma_a\partial_t {\bf P}_t^a +
\omega_a^2 {\bf P}_t^a = \omega_{pa}^2 {\bf E}_t\ , \ \ \ \ \ 
{\bf P}_t = \sum_{a=1}^N {\bf P}_t^a\ .
\ee
Together with Maxwell's equations, 
Eq. (\ref{p}) form a system of sought-for causal evolution equations which are
to be transformed into a system of first order 
equations by means of the Hamiltonian formalism. In finite difference
time domain numerical schemes, the Hamiltonian formalism
for the Lorentz model has been used in \cite{hflm} to study
propagation of an electromagnetic pulse in homogeneous Lorentz media.

In our approach no boundary conditions are imposed on electromagnetic
fields at medium and/or target interfaces. 
The latter are modeled by spatially
dependent couplings of media to electromagnetic fields which are
included into the system Hamiltonian. At any
interface, the couplings are allowed to have discontinuities, or,
from a physical point of view, they remain smooth but change
rapidly, $\lambda_w|\mbf{\nabla }\omega_p|/\omega_p >\!\! > 1$,
at the interface, where $\lambda_w$ is a typical wave length
of the incoming wave packet. The conventional boundary conditions
are automatically generated by the dynamics \cite{landau}. 
Thus, the initial
value problem is solved in ${\mathbb L}_2({\mathbb R}^3)$ for every
matter and electromagnetic field component. This
implies that the energy of the propagating wave packet remains
finite (in contrast to the scattering matrix approach based on 
plane wave solutions).
 
Let us now formulate the initial value problem for
a generic passive medium and then apply the formalism
to multi-resonant Lorentz models. 
Combine the fields, inductions and medium responses into columns:
$$
\psi_t^F = \pmatrix{{\bf E}_t \cr {\bf H}_t}\ ,\ \ \ \ 
\psi_t^I = \pmatrix{{\bf D}_t \cr {\bf B}_t}\ ,\ \ \ \ 
\psi_t^R = \pmatrix{{\bf P}_t \cr {\bf M}_t}\ .\ \ \ \ 
$$
Assuming linear response theory, one can write for the Fourier transforms 
\be
\label{response}
\tilde{\psi}^R_\omega = \tilde{\chi}_\omega\tilde{\psi}_\omega^F\ ,
\ee
where the Fourier transform of a general response function,
$\tilde{\chi}_\omega$, has to satisfy a dispersion relation that ensures causality (like the 
Kramer-Kronig relations for the dielectric constant) \cite{landau}. 
For anisotropic media, $\tilde{\chi}_\omega$ is a symmetric matrix
acting on components of electromagnetic fields.
With this type of generality all possible media are
covered as long as linear response theory is valid. 
The response function $\tilde{\chi}_\omega$ can either be
modeled or
measured and tabulated in some frequency range of interest (determined
by the frequency bandwidth of the initial wavepacket), say, $\omega
\in [\omega_1,\omega_2]$. Next, the components of
$\tilde{\chi}^{-1}_\omega$ are expanded in a basis
of suitable orthogonal polynomials. An optimal expansion is often achieved
in the Chebyshev polynomial basis. Chebyshev polynomials are defined
in the interval $[-1,1]$ so a corresponding rescaling and
translation of $[\omega_1, \omega_2]$ must be done.
By taking the Fourier transform of $\tilde{\chi}^{-1}_\omega \tilde{\psi}^R_\omega = 
\tilde{\psi}_\omega^F$ we obtain the desired differential equation
\be
\label{gm}
\sum_{n=0}^N \chi_n \partial_t^n \psi^R_t = \omega_p \psi_t^F\ ,
\ee
where $\omega_p=\omega_p({\bf r})$ plays the role of the coupling constant between
matter and electromagnetic fields.
The order $N$ is determined by the  
highest order of polynomials used to approximate $\tilde{\chi}_\omega^{-1}$.
The expansion coefficients $\chi_n$ and the coupling $\omega_p$ are matrices for
anisotropic media.

The basic idea of the Hamiltonian formalism is to convert the 
system (\ref{p}) or (\ref{gm}) into a system
of first-order differential equations by introducing auxiliary
(matter) fields. 
The number of such fields is determined by
the order of the original evolution equation for matter. 
For instance, in the case of the multi-resonant Lorentz model, 
there are
$N$ fields ${\bf P}_t^a$ each of which satisfies a second order 
differential equation. In the Hamiltonian
formalism one would have $2N$ real vector fields, $\mbf{\xi}_t^j$, 
$j=1,2,...,2N$. A simple possibility is to set
\ba
\label{aux}
{\bf P}_t^a &=&(\omega_{pa}/\omega_a)\, \mbf{\xi}_t^{2a-1} \ ,\\
\label{aux1}
\partial_t\mbf{\xi}_t^{2a-1} &=& \omega_a\mbf{\xi}_t^{2a}\ ,\\
\label{aux2}
\partial_t\mbf{\xi}_t^{2a}  
&=& -2\gamma_a\mbf{\xi}_t^{2a} -\omega_a\mbf{\xi}_t^{2a-1} + \omega_{pa}{\bf E}_t\ .
\ea
The reason of inserting the factor $\omega_{pa}/\omega_a$ in the definition
(\ref{aux}) of 
the auxiliary fields  will be evident from
what follows. Note that the medium polarization ${\bf P}_t$ must be zero
in empty space where $\omega_{pa} = 0$. The factor $\omega_a^{-1}$
in (\ref{aux}) simplifies the energy conservation and stability
analysis.
 
For the Lorentz model there is another convenient way to introduce 
the Hamiltonian formalism by using $N$ complex
vector fields $\mbf{\zeta}_t^a$ which satisfy the first order differential equation
\ba
\label{caux}
\partial_t\mbf{\zeta}_t^a &=& \lambda_a \mbf{\zeta}_t^a - i\omega_{pa}{\bf E}_t\ ,\\
\label{caux2}
{\bf P}_t^a &=& \frac{\omega_{pa}}{2\nu_a}\ \left(\mbf{\zeta}_t^a + \bar{\mbf{\zeta}}_t^a\right)\ ,
\ea
where $\lambda_a = -\gamma_a + i\nu_a$ and $\nu_a = \sqrt{\omega_a^2 -\gamma_a^2}$. 
This representation is defined only if $\gamma_a < \omega_a$ (i.e., the attenuation is
not high).
From the numerical
point of view, solving a {\it decoupled} system of 
$N$ first order differential equation and taking complex conjugation 
(denoted here by an over bar)  is less expensive 
than solving an original system of differential equations
to compute the medium polarization. 

Returning to the  general case, we introduce a set of auxiliary fields
$\xi_t$
to convert (\ref{gm}) into a first-order system,
\be
\label{matter}
\partial_t \xi_t = {\cal H}_M^{F} \xi_t + {\cal V}_{MF} \psi_t^F\ .
\ee
The operators ${\cal H}_M^{F}$ and ${\cal V}_{MF}$  are 
determined by the details of going over to the Hamiltonian formalism.
We shall call ${\cal H}_M^{F}$ the matter Hamiltonian; 
it governs time evolution of the medium when no external fields
are applied. The index $F$ indicates that the electromagnetic 
degrees of freedom are described by fields, not
inductions. We shall see shortly that the matter Hamiltonian 
depends on whether $\psi_t^I$ or $\psi_t^F$ is used
as independent electromagnetic variables. 
The matrix ${\cal V}_{MF}$ describes the coupling of matter 
to the electromagnetic fields, which is emphasized
by the index $MF$ (matter-to-field coupling). 
We introduce a linear time independent 
operator ${\cal R}$ that acts in the space of auxiliary (matter) fields so that
\be
\label{respoper}
\psi_t^R = {\cal R} \xi_t\ ,
\ee
that is, the (response) operator ${\cal R}$ maps a 
given configuration of auxiliary fields onto the 
corresponding physical response
field. It depends on the definition of the matter fields
(cf. (\ref{aux}) and (\ref{caux2})).
A passive medium is not excited, $\xi_t=0$, if no external
electromagnetic field is applied; that is, the initial
condition for  Eq. (\ref{matter}) is such that it has only
the trivial solution whenever $\psi_t^F=0$. Under this
condition, the solution of (\ref{matter}) reads
$$
\xi_t = \int_{-\infty}^t d\tau e^{(t-\tau){\cal H}_M^F}{\cal V}_{MF} 
\psi_\tau^F\ .
$$
Hence, the linear response operator $\tilde{\chi}_\omega$ in 
(\ref{response}) is the Fourier transform of the operator
\be
\label{resp11}
\chi_t=\theta_t{\cal R}\,e^{t{\cal H}_M^F}\,{\cal V}_{MF}\ ,
\ee
where $\theta_t$
is the Heaviside function. Or, vice versa,  ${\cal R},\, {\cal H}_M^F$
and ${\cal V}_{MF}$ must be chosen so that the Fourier transform
of $\chi_t$ defined by (\ref{resp11}) coincides with the known response function
$\tilde{\chi}_\omega$ of the medium in a designated frequency range.
 
Maxwell's equations without external currents can be rewritten in the Hamiltonian form
\be
\label{field}
\partial_t\psi_t^F = 
{\cal H}_F\psi_t^F -\partial_t \psi_t^R = {\cal H}_F\psi_t^F + {\cal V}_{FM} \xi_t\ .
\ee
The field-to-matter coupling ${\cal V}_{FM}$ and the field Hamiltonian ${\cal H}_F$ are 
deduced from (\ref{matter}) by acting on the latter by the operator
${\cal R}$, which yields
\ba
\label{vfm}
{\cal V}_{FM} &=& - {\cal R}{\cal H}_M^{F}\ ,\\
\label{hf}
{\cal H}_F &=& \pmatrix{0 & c\mbf{\nabla}\times \cr -c\mbf{\nabla}
\times & 0} 
- {\cal R}{\cal V}_{MF}\equiv {\cal H}_0 -{\cal R}{\cal V}_{MF}\ .
\ea
It is always possible to set up the Hamiltonian formalism 
so that ${\cal R}{\cal V}_{MF} \equiv 0$ and, hence, ${\cal H}_F={\cal H}_0$.
It is not difficult to verify that this holds for the Lorentz model discussed above. In the general case,
the standard procedure of going over to the Hamiltonian formalism \cite{arnold},
where components of $\xi_t$ are identified with time derivatives of
the response field, 
$\xi_t^k \sim \partial_t^k \psi_t^R$,
leads to the same result that ${\cal R}{\cal V}_{MF}=0$. Thus, without loss of generality, the last
term in the field Hamiltonian (\ref{hf}) can be omitted.

The auxiliary matter and electromagnetic fields (or inductions) are unified into a larger column
\be
\label{state}
\Psi_t^F = \pmatrix{ \psi_t^F\cr \xi_t}\ ,\ \ \ \ \Psi_t^I =\pmatrix{\psi_t^I\cr \xi_t}\ .
\ee
The wave function $\Psi_t^F$ satisfies the Schr\"odinger equation
\be
\label{main}
\partial_t \Psi_t^F = {\cal H}^{F} \Psi_t^F\ , \ \ \ \ 
{\cal H}^{F} = \pmatrix{{\cal H}_0 & {\cal V}_{FM}\cr {\cal V}_{MF} &
  {\cal H}_M^{F}}\ .
\ee
which has to be solved with the initial field configuration 
$\psi_{t=0}^F = \psi_0$, while the matter fields
are assumed to be zero at the initial moment of time, $\xi_{t=0} = 0$, e.g., 
the initial wave packet is localized in an empty space region.
Equations  (\ref{matter})
and (\ref{field}) are equivalent to (\ref{main}). In a similar fashion, one can derive 
the Schr\"odinger equation
for $\Psi_t^I$. Note that
\be
\label{ftoi}
\Psi_t^I = {\cal S}\Psi_t^F\ ,\ \ \ \ \ 
{\cal S} = \pmatrix{1 & {\cal R} \cr 0 & 1 }\ ,\ \ \ \ 
 {\cal S}^{-1} = \pmatrix{1 & -{\cal R} \cr 0 & 1 }\ ,
\ee
Hence,
\be
\label{main2}
\partial_t \Psi_t^I = {\cal H}^{I} \Psi_t^I\ ,\ \ \ \ \ 
{\cal H}^{I} = {\cal S}{\cal H}^{F}{\cal S}^{-1}\ .
\ee
The corresponding blocks of ${\cal H}^{I}$ have the form
\ba
\label{IF1}
{\cal H}_I &=& {\cal H}_0 \ ,\ \ \ \ \ \ \ \ \ {\cal V}_{MI} = {\cal V}_{MF}\ ,\\
\label{VIM}
{\cal V}_{IM} &=& {\cal V}_{FM} + 
{\cal R}{\cal H}_M^F - {\cal H}_0{\cal R} = -{\cal H}_0{\cal R}\ ,\\ \label{HI}
{\cal H}^{I}_M &=& {\cal H}^{F}_M - {\cal V}_{MF}{\cal R}\ .
\ea
To simplify ${\cal V}_{IM}$, Eq. (\ref{vfm}) has been applied. 
Observe in (\ref{HI})  the aforementioned dependence 
of the matter Hamiltonian on the representation of electromagnetic degrees of freedom.
The use of either (\ref{main}) or 
(\ref{main2}) in numerical simulations has its own advantages and
disadvantages 
which are discussed below.   

As an example, we give an explicit form of the Hamiltonian for the Lorentz model when the auxiliary
field are defined by (\ref{aux})
\ba
\label{lmvfm}
{\cal V}_{FM} &=& ({\cal V}_{FM1}, {\cal V}_{FM2}, \cdots, {\cal V}_{FMN})\ ,\ \ \ \ {\cal V}_{FMa} = 
\pmatrix{0 & -\omega_{pa}  \cr 0&0}\ ,\\ \label{lmvmf}
{\cal V}_{MF} &=&- {\cal V}_{FM}^*\ ,\\ \label{lmhm}
{\cal H}_M^{F} &=&{\sf diag\ }\left({\cal H}_{M1}^{F},\ {\cal
H}_{M2}^{F}, 
\cdots, {\cal H}_{MN}^{F}\right)\ , \ \ \ \ 
{\cal H}_{Ma}^{F} = \pmatrix{0 & \omega_a \cr -\omega_a & -2\gamma_a}\ ,
\ea
where ${\sf diag}$ indicates that the corresponding matrix 
is block-diagonal with blocks listed in the order 
from the upper left to lower right corners. Note that the matrices
${\cal V}_{FMa}$ and ${\cal H}_{FMa}$ act on a six-dimensional column
$(\mbf{\xi}^{2a-1}_t,\mbf{\xi}^{2a}_t)^*$. Therefore they should be 
understood as composed of $3\times 3$ blocks. Each block is obtained
by multiplying the unit matrix by the number indicated in place of the
block in (\ref{lmvfm}) and (\ref{lmhm}).

Our final remark in this section concerns ``canonical'' 
transformations in the Hamiltonian formalism. As has been pointed
out, the auxiliary fields are not uniquely defined. There is a freedom
of making general complex nonsingular 
linear transformations such as
\be
\label{ct}
\xi_t \rightarrow {\cal S}_M \xi_t\ ,\ \ \ \ \det{\cal S}_ M \neq 0\ .
\ee
If the infinitesimal evolution operator ${\cal U}_{\Delta t}^{(F,I)}=
\exp(\Delta t {\cal H}^{(F,I)})$ 
is computed with one choice of the auxiliary fields, a simple
similarity transformation, 
like the one in (\ref{main2}), would allow us to compute it in any other 
basis of auxiliary fields. This is an important observation because
the auxiliary field 
basis can be chosen  in a way that
facilitates computation of the evolution operator (e.g., to improve the
convergence rate or speed up simulations). For instance,
in the complex representation (\ref{caux2}) of the auxiliary 
fields in the Lorentz model, the matter Hamiltonian is
diagonal. The corresponding transformation of the auxiliary fields is given by
\be
\label{S}
\pmatrix{\mbf{\xi}^{2a-1}_t\cr \mbf{\xi}^{2a}_t} = \frac{1}{2\nu_a}\, 
\pmatrix{\omega_a& \omega_a\cr \lambda_a &\bar{\lambda}_a}
\pmatrix{\mbf{\zeta}_t^a \cr \bar{\mbf{\zeta}}^a_t} 
\equiv {\cal S}_M \pmatrix{\mbf{\zeta}_t^a \cr \bar{\mbf{\zeta}}^a_t}\ .
\ee
To transform the whole system into this representation, the Hamiltonian ${\cal H}^{F}$ is 
replaced by ${\cal S}^{-1}{\cal H}^{F}{\cal S}$ and the wave function $\Psi_t^F$ by
${\cal S}\Psi_t^F$ where ${\cal S}$ is block-diagonal with the unit matrix in the upper left (field) corner
and with ${\cal S}_M$ in the lower right (matter) corner.

\section{The grid representation of Maxwell's theory}
\setcounter{equation}0

Consider an equidistantly spaced finite grid with periodic boundary conditions. Let $\Delta r$ be
the grid step and ${\bf n}$ be a vector with integer valued
components. Then the dynamical variables are projected onto the grid
by taking their values at grid points ${\bf r}= {\bf n }\Delta r$,
\be
\label{grid1}
\Psi_t^Q({\bf r}) \rightarrow \Psi_t^Q({\bf n}\Delta r)\ ,
\ee
where $Q$ denotes the representation, $I$ or $F$. For simplicity, 
a cubic grid is assumed here. It is straightforward to generalize
the discussion to a generic rectangular grid.
Consider a discrete
Fourier transformation associated with the grid \cite{fft,lill}
\be
\label{grid2}
\tilde{\Psi}_t^Q ({\bf n}k_0) = \sum_{{\bf n}^\prime}{\cal F}_{{\bf
nn}^\prime} \Psi_t^Q({\bf n}^\prime\Delta r)\ ,\ \ \ \ {\cal F}^*{\cal F}= {\cal
F}{\cal F}^* = 1\ ,
\ee
where the dual lattice step is $k_0= 2\pi/\Delta r$. The grid spatial size
$L$ and step must be chosen so that the Fourier transform of the initial
wavepacket has support within the region
$k\in [k_{min}, k_{max}]$ where $k=|{\bf k}|$,  $k_{max} = k_0$ and $k_{min} = 2\pi/L$.
The Hamiltonian ${\cal H}^Q$ is split into a sum 
\be
\label{grid3}
{\cal H}^Q = {\cal H}_0^Q + {\cal V}^Q\ ,
\ee
where all the spatial derivatives are included into ${\cal H}_0^Q$ and
${\cal V}^Q$ contains multiplications by position dependent functions.
This is always possible for the Hamiltonian described in the
preceding section. The operator ${\cal V}^Q$ is projected naturally
\be
\label{grid4}
{\cal V}^Q({\bf r})\Psi_t^Q({\bf r})
\rightarrow {\cal V}^Q({\bf n}\Delta r) \Psi_t^Q({\bf n}\Delta r)\ .
\ee
Consider ${\cal H}^Q_0$ in the Fourier basis, ${\cal
H}^Q_0(\mbf{\nabla})\rightarrow {\cal H}_0^Q(i{\bf k})$. The
projection is then done via the discrete Fourier transform
\be
\label{grid5}
\left.
{\cal H}_0^Q(\mbf{\nabla}) \Psi_t^Q({\bf r})\right|_{{\bf r}= {\bf
n}\Delta r}
\rightarrow\sum_{{\bf n}^\prime} \left({\cal F}^*\right)_{{\bf
n}{\bf n}^\prime} {\cal H}_0^Q(i{\bf n}^\prime k_0)
\tilde{\Psi}_t^Q({\bf n}^\prime k_0)\ .
\ee
In what follows, the rules (\ref{grid4})
and (\ref{grid5}) {\it define} the action of the operators
${\cal V}^Q$ and ${\cal H}^Q_0$ and their functions on any state
vector. The action of a product of ${\cal V}^Q$ and ${\cal H}^Q_0$
on any state vector is understood as consecutive actions of
these operators according to the rules (\ref{grid4}) and (\ref{grid5}),
in the order specified in the product.

The projection (\ref{grid5}) as well as any action of ${\cal H}_0^Q$
on state vectors 
is performed by the fast Fourier method \cite{fft}. It requires $N\log_2 N$
elementary operations (flops)
with $N$ being the grid size. 
In finite differencing schemes,
the action of ${\cal H}_0^Q$ on a state vector 
would require $mN_d$ operations where the integer
$m$ depends on a particular difference scheme used to approximate
derivatives, and $N_d$ is the grid size used in the differencing
scheme. It should be noted that, as shown below, the use
of the fast Fourier transform eliminates the phase error (because
the correct electromagnetic dispersion relation is preserved)
and operates at the Nyquist limit. These two features allows one
to reduce substantially the grid size as compared with that 
in a finite differencing scheme, while providing the same accuracy
in simulations. Recall that, in scattering problems, the phase of the return
signal contains the most significant information about the target.
So, in practice, grids in finite differencing schemes are significantly
larger (more dense) than grids used in the fast Fourier method in order to
reduce the phase errors in the former. Needless to say, the advantage of the fast 
Fourier method in reducing the phase error becomes even more 
significant in higher dimensions because
$N_d/N = (n_d/n)^D$ where $n_d$ and $n$ are the corresponding numbers of grid points
per shortest wave length in the wave packet, and $D$ is the grid dimension.
The Nyquist limit is $n=2$, while $n_d$ is of order 10 or higher.

\section{The split operator method}
\setcounter{equation}0

Let $|\Psi|$ denote the ${\mathbb L}_2({\mathbb R}^3)$ norm of the wave function, or the
Euclidean norm of the corresponding vector (\ref{grid1}) 
in the grid representation.
One possible way to compute numerically the path integral (\ref{pi}) 
is based on the Kato-Trotter product formula \cite{r}
\be
\label{kato}
\lim_{n\rightarrow\infty} \left|e^{t({\cal A} +{\cal B})}\Psi - 
\left(e^{t{\cal A}/2n}e^{t{\cal B}/n}e^{t{\cal A}/2n}\right)^n\Psi\right| =
0\ ,
\ee
for a general $\Psi$ and under certain assumptions about 
the linear operators ${\cal A}$ and ${\cal B}$ in the Hilbert space 
spanned by $\Psi$. For our purposes it is sufficient to note that
for bounded operators, (\ref{kato}) always holds and is known as the Lie-Trotter
product formula. In the grid representation, which would always be assumed,
unless stated otherwise,  operators ${\cal A}$ and
${\cal B}$ are finite matrices and, hence, bounded.

Let us apply (\ref{kato}) to the split (\ref{grid3}), meaning that
the operator
${\cal H}_0^{Q}$ is used in place of ${\cal A}$ (or ${\cal B}$) and, respectively,
the operator
${\cal V}^{Q}$ is used in place of ${\cal B}$ (or ${\cal A}$).
The infinitesimal evolution operator in (\ref{te}) can be approximated by the first
term in the following expansion
\ba
\label{split}
{\cal U}_{\Delta t}^Q&=&
e^{\Delta t({\cal A} +{\cal B})} = {\cal G}_{\Delta t}^Q + 
\Delta t^3{\cal W}_{\Delta t}\\
\label{splitg}
{\cal G}_{\Delta t}^Q &=&
e^{\Delta t{\cal A }/2}e^{\Delta t{\cal B}}e^{\Delta t{\cal A }/2} \\
\label{splitw}
{\cal W}_{\Delta t}&=&
-\frac{1}{24}\left( [{\cal A},[{\cal A},{\cal B}]] -2
[{\cal B},[{\cal B},{\cal A}]]\right) + O(\Delta t)\ .
\ea
By making 
$n$ larger while keeping $n\Delta t = t$ fixed, the strong convergence in
(\ref{kato}) guarantees that the error can be made arbitrary small for
any initial state, 
\be
\label{splitconv}
|({\cal U}_{\Delta t}^Q)^n \Psi_0 - ({\cal G}_{\Delta
t}^Q)^n\Psi_0|\rightarrow 0
\ee
as $n\rightarrow \infty$ for any $\Psi_0$.
The numerical iteration algorithm is then based on the replacement of the exact
evolution (\ref{te}) by the approximate one
\be
\label{teapp}
\Psi_{t+\Delta t}^Q = {\cal G}_{\Delta t}^Q \Psi_t^Q\ .
\ee  
The quantity $t\Delta t^2|{\cal W}_0^Q\Psi_0|/|\Psi_0|$ 
can be used to roughly estimate
the accuracy of the algorithm.  A more detailed accuracy analysis is
given in Section 8.
By making use of the Campbell-Hausdorf
formula for the exponential of the sum of operators 
it is possible to obtain the symmetric product formula in
(\ref{split}) to approximate ${\cal U}_{\Delta t}$
up to any desired order in $\Delta t$, presumably achieving a higher 
accuracy \cite{raedt}. 
This would come at the price of having more
exponentials in the symmetric product ${\cal G}_{\Delta t}^Q$.
In numerical simulations, 
one should keep in mind that  computational costs of decreasing
$\Delta t$ in the third order split (\ref{split}) (i.e., increasing the 
number of steps in the time evolution) might be less 
than those of computing a lesser number of actions of ${\cal
G}_{\Delta t}^Q$ in higher order splits. So, the higher order splits
are not always optimal to achieve a better accuracy \cite{kosloff}.

On the grid, the action of the amplification operator
${\cal G}_{\Delta t}^Q$ is computed according to the rules
(\ref{grid4}) and (\ref{grid5}) applied to, respectively,
$\exp(\Delta t {\cal V}^Q)$ and $\exp(\Delta t {\cal H}^Q_0)$. 
Explicit formulas for the exponentials of the corresponding operators can
be worked out in the field and induction representations.
If the fields are used as independent variables, then a natural choice is  
\be
\label{fsplit}
{\cal H}^{F} = {\cal H}_0^{F} + {\cal V}^{F} = \pmatrix{{\cal H}_0 &0 \cr 0&{\cal H}_M^{F}} +
\pmatrix { 0 & {\cal V}_{FM}\cr {\cal V}_{MF} & 0}\ .
\ee 
The matter Hamiltonian ${\cal H}_M^F$ can also be transferred into
${\cal V}^{F}$ 
if so desired. This rearrangement
affects the accuracy of the method, meaning that the operator (\ref{splitw})
would change. In turn, a rearrangement of operators in the split can be used to improve the accuracy.
We shall discussed this issue later. Using 
the Taylor series we infer that
\ba
\label{freeHF}
\exp(t{\cal H}_0^{F})& =& \pmatrix{ \exp({t{\cal H}_0}) & 0 \cr 0 & \exp({t{\cal H}_M^{F}})}\ ,\\
\label{freeH}
\exp(t{\cal H}_0) &=& 1 + \left[\cos(ct\sqrt{-\Delta}) - 1\right]{\cal
P}_\perp + 
\frac{\sin\left(ct\sqrt{-\Delta}\right)}{c\sqrt{-\Delta}}\ {\cal H}_F\ ,
\ea
where $\Delta = \mathbf{\mbf{\nabla}}\cdot\mathbf{\mbf{\nabla}}$ is the Laplace operator and 
${\cal P}_\perp = 1 - {\mathbf{\mbf{\nabla}}} (\Delta^{-1}\mbf{\nabla}
\cdot\ )$ 
is the projector on transverse fields, that is,
${\cal P}_\perp {\bf E} = {\bf E}$ if 
$ \mbf{\nabla}\cdot {\bf E} = 0$ and $P_\perp {\bf E} = 0$ if the vector field
${\bf E}$ is conservative, ${\bf E} = \mbf{\nabla} \phi$. The
projector ${\cal P}_\perp$ can be 
omitted in (\ref{freeH}) if it is known (e.g. from a theoretical analysis of the system)
that the fields remain transversal in due course. In this case, the two first terms in
(\ref{freeH}) are equal to $\cos(ct\sqrt{-\Delta})$.
The action of  $\exp(t{\cal H}_0) $ is computed 
by the fast Fourier transform according to (\ref{grid5}). 
In the Fourier basis, $-\Delta = {\bf k}^2= {\bf n}^2 k_0^2$. 
Note also that the Fourier
transform of the fields $\psi_t^F$ is required, 
while the auxiliary fields remain in the grid basis all the time.
The exponentials of ${\cal H}^{F}_M$ and ${\cal V}^{F}$ can either be computed analytically for
simple models like the Lorentz model, as is shown Section 5, or, in general 
case, by direct diagonalization at each grid site. 

Alternatively, 
the following approximation of the exponential of an operator can be used
\be
\label{inv}
e^{\Delta t {\cal B}} = \frac{1 +\Delta t {\cal B}/2 +
\Delta t^2 {\cal B}^2/12}{1 -\Delta t {\cal B}/2 +
\Delta t^2 {\cal B}^2/12} +O(\Delta t^5) =
\frac{1 +\Delta t {\cal B}/2}
{1 -\Delta t {\cal B}/2} +O(\Delta t^3) \ .
\ee
If the matrix ${\cal B}$ is anti-hermitian, then the approximations
(\ref{inv}) of the exponential of ${\cal B}$ retain unitarity, which
is important for stability of the split algorithm (see Theorems 7.1
and 7.3). On the other hand, costs of computing the inverse
matrices
in the right hand side of (\ref{inv}) can be less than those of 
computing the exponential.

If the inductions are used as independent variables, then a natural choice of the split would be
 \be
\label{isplit}
{\cal H}^{I} = {\cal H}_0^{I} + {\cal V}^{I} 
= \pmatrix{{\cal H}_0 & - {\cal H}_0{\cal R} \cr 0& 0} +
\pmatrix { 0 & 0 \cr {\cal V}_{MI} & {\cal H}_M^{I}}\ .
\ee 
 Making use of the Taylor expansion again we deduce that
\be
\label{eh0i}
\exp(t{\cal H}_0^{I}) = \pmatrix{ \exp({t{\cal H}_0}) 
& \left[1-\exp(t{\cal H}_F)\right]{\cal R}\cr 0& 1}\ .
\ee
and, similarly,
\be
\label{eVI}
\exp(t{\cal V}^{I})= \pmatrix{1 & 0\cr ({\cal H}_M^{I})^{-1}
[\exp(t{\cal H}_M^{I}) - 1]{\cal V}_{MI} &
\exp(t{\cal H}_M^{I})}\ .
\ee

Now we can compare the two splits.  The split (\ref{fsplit})
has an advantage over (\ref{isplit}) because it requires less calls of the fast Fourier transform. 
Indeed, in the former the fast Fourier transform is called only for
the fields $\psi_t^F$. As it follows from (\ref{eh0i}), 
the operator $\exp(t{\cal H}_0)$ acts on both the
inductions $\psi_t^I$ and the auxiliary fields. 
Hence the fast Fourier transform must be called for the entire column $\Psi_t^I$.
If the number of auxiliary fields is large, there 
might be a substantial difference in the computational speed 
of two algorithms. The latter, however, depends on the choice of
the matter fields which, in turn, determines ${\cal R}$ and therefore
the number of calls of the fast Fourier transform.
Note that the matter fields can always be chosen in such a way that 
only one of the components of $\xi_t$ specifies the response field
$\psi_t^R$. Thus, the canonical transformation (\ref{ct}) can be used
to reduce the number of calls of the fast Fourier transform.
If ${\cal R}$ is chosen so that it depends on position, multiplication of the
matter fields by ${\cal R}$ must be done before calling the fast
Fourier transform. 
A significant advantage of the split in the induction representation is that the
Gauss law can be exactly fulfilled without altering the algorithm
(see Theorem 8.2).

In empty space either of the splits reproduces an {\it exact} solution
of Maxwell's equations for any period of time $t$, provided the
initial pulse is bandwidth limited. Indeed, on the grid, the initial
wave packet is a superposition of a finite number of plane waves.
Thanks to the linearity of the theory, each Fourier
mode is evolved exactly, that is, {\it without} any phase error, by
$\exp(t{\cal H}_0)$ for any $t>0$. 
As final remarks in this section, we note that the algorithm can operate
at the Nyquist limit: Two grid points per shortest wavelength 
in the initial wave packet \cite{fft}. Yet, for the multiresonant
Lorentz model, it is {\it unconditionally} stable 
(see Theorems 7.1 and 7.3). These features cannot be 
achieved in any finite difference scheme.

\section{A multi-resonant Lorentz model}

An analytical expression for the exponents of
 ${\cal H}_M^{Q}$ and ${\cal V}^{Q}$  helps to reduce computational costs.
Here such analytical expressions are derived for multi-resonant 
Lorentz models.
Let us take first the field representation.
Due to the block diagonal structure of ${\cal H}_M^F$ we get
\ba
\label{eHM2}
\exp(t{\cal H}_M^{F}) &=& {\sf diag\ }\left( \exp(t{\cal H}_{M1}^{F}),\ \exp(t{\cal H}_{M2}^{F}), 
\ \cdots,\  \exp(t{\cal H}_{MN}^{F})\right)\ ,\\
\label{eHM}
\exp(t{\cal H}_{Ma}^{F}) &=& e^{-\gamma_at}\left[\cosh \tilde{\nu}_a t +
\frac{\sinh\tilde{\nu}_at}{\tilde{\nu}_a}\left({\cal H}_{Ma}^{F} +\gamma_a
\right)\right]\ ,
\ea
where $\tilde{\nu}_a = (\gamma_a^2 - \omega_a^2)^{1/2}$.
The exponential (\ref{eHM}) is easy to compute 
by expanding ${\cal H}_{Ma}^{F}$ in the Pauli matrix basis, which is
also a basis for the Lie algebra $su(2)$,
and then by using the well known formula for the exponential of a linear combination of Pauli matrices.
For small attenuation, $\gamma_a < \omega_a$, we get $\tilde{\nu}_a =
i\nu_a$. The hyperbolic
functions in (\ref{eHM}) become trigonometric ones and
$\tilde{\nu}_a$ is replaced by $\nu_a$. 
The eigenvalues of the matter Hamiltonian are $\lambda_a = -\gamma_a
\pm \tilde{\nu}_ a$. Hence, ${\rm
Re}\, \lambda_a < 0$ and amplitudes of the matter fields are always
exponentially attenuated as $t\rightarrow \infty$, unless $\gamma_a
=0$ leading to ${\rm Re}\, \lambda_a =0$.

Computation of $\exp(t{\cal V}^{F})$ is a bit more subtle. We start with the observation that
the characteristic polynomial of ${\cal V}^{F}$ has a simple form
\be
\label{omegap}
\det \left({\cal V}^{F} - \lambda\right) = \lambda^{2N}(\lambda^2 +
\omega_p^2)\ ,\
\ \ \ \ \ 
\omega_p^2 = \sum_{a=1}^N \ \omega_{pa}^2\ .
\ee
This can be proved either by a direct computation or by mathematical induction.
So, ${\cal V}^F$ has  $2N$ zero eigenvalues and two non-zero ones, $\lambda =\pm i\omega_p$.
Let $X$ be the eigenvector of ${\cal V}^F$ corresponding to the eigenvalue $i\omega_p$.
Its components have the form
$$
X_j = \omega_p^{-1}\left[ {\cal V}^{F}_{j1} + i\omega_p
\delta_{j1}\right]\ ,\ \ \ j=1,2,..., 2(N+1)\ ,
$$
so that  $\bar{X}\cdot X = 1$ and $\bar{X}\cdot\bar{X} = X\cdot X =0$
where the dot denotes the Euclidean scalar product. 
The skew-symmetric matrix ${\cal V}^{F}$ has the following spectral 
decomposition
\be
\label{v}
{\cal V}^{F} = i\omega_p \left( X\otimes \bar{X} - \bar{X}\otimes X\right)\ .
\ee
Taking the square of (\ref{v}) we also infer that
\be
\label{XX}
X\otimes \bar{X} = -\omega_p^{-2} \left({\cal V}^{F}\right)^2 -i\omega_p^{-1}{\cal V}^{F}\ .
\ee
The exponential of (\ref{v})  is obtained via the Taylor series
and making use of (\ref{XX}).
The final result reads
\be
\label{etvf}
\exp(t{\cal V}^{F}) = 1 +\frac{\sin\omega_pt}{\omega_p}\ {\cal V}^{F}
+ 2\left(\frac{\sin(\omega_pt/2)}{\omega_p}\right)^2 \, \left({\cal V}^{F}\right)^2 \ .
\ee

In the induction representation, an explicit formula for $\exp(t{\cal H}_M^{I})$ is not
that simple. To avoid unnecessary technicalities,
we limit the discussion to the simplest case of the one-resonant
Lorentz model. We choose the matter fields so that $\mbf{\xi}_t^1 =
{\bf P}_t$ and $\mbf{\xi}_t^2 = \partial_t \mbf{\xi}_t^1$. In this
case, non-zero elements of the matter Hamiltonian are ${\cal
H}^I_{M12} = 1$, ${\cal H}^I_{M21} = -\omega_0^2 -\omega^2_p$, and
${\cal H}^I_{M22}= -2\gamma$. The coupling matrix ${\cal V}_{MI}$
has only one non-zero element, ${\cal V}_{MI21}= \omega_p^2$.
Here $\omega_0$ is the resonant frequency, $\gamma$ is the
attenuation constant and $\omega_p$ is the plasma frequency.
Using the Pauli matrix basis again, we find that
the expression (\ref{eHM}) holds for $\exp(t{\cal H}_M^{I})$
if we replace in it $\nu_a$ by $\nu_{p} = \sqrt{\omega_0^2
+\omega_{p}^2 -\gamma^2 }$, $\gamma_a$ by $\gamma$ and ${\cal H}_M^F$ 
by ${\cal H}_M^I$.
 The lower left corner of (\ref{eVI}) has the form
\be
\nonumber
\left({\cal H}_M^{I}\right)^{-1}\left(\exp(t{\cal H}_{M}^{I}) - 1\right) {\cal V}_{MI}
=-\frac{\omega_{p}}{\omega_0^2 + \omega_{p}^2}\left(\exp(t{\cal H}_{M}^{I}) - 1\right)\
\pmatrix{1&0\cr 0&0}\ .
\ee
A further simplification can be achieved by going over to the complex representation (\ref{caux})
of the auxiliary fields in which the matter Hamiltonian is diagonal. 
The transformation rule is explained
in the paragraph after Eq. (\ref{S}). 

\section{Energy and norm conservation}
\setcounter{equation}0

Consider the ${\mathbb L}_2({\mathbb R}^3)$ 
norm of $\Psi_t^Q$,  $|\Psi_t^Q|^2 = \int d{\bf r}
\Psi_t^{Q*}\Psi_t^Q \equiv (\Psi_t^Q,\Psi_t^Q)$. 
In the grid representation, 
the norm coincides with the corresponding (complex) Euclidean norm,
$(\Psi^Q_t,\Psi_t^Q) = 
\Delta r $ $\sum_{\bf n} \Psi_t^{Q*}({\bf n}\Delta r)\Psi_t^Q({\bf
n}\Delta r)$ where the sum is taken over all grid sites. 
By taking the time
derivative of $|\Psi_t^Q|^2$
and using the evolution equation (\ref{sch}), it is not hard to deduce
that the norm is conserved, provided the Hamiltonian is anti-Hermitian
\be
\nonumber
{\cal H}^{Q*} = - {\cal H}^Q\ .
\ee
For a generic passive media this is not the case. So the norm 
is generally not
conserved in contrast to the quantum mechanical case. However, we
shall see that in the case when the matter evolution is described
by second order differential equations in time and no attenuation
is present, the norm coincides with the system energy and is
conserved. In numerical simulations, this important property
can be used to help to control the accuracy.

Consider
multi-resonant Lorentz models with no attenuation, $\gamma_a=0$.
We start with the observation that the field and matter evolution
equations can be obtained from the variational principle for the action
\be
\label{action}
S = \int dt L = \int\! dt\!\!\int\! d{\bf r}\left[
\frac 12 \left( {\bf E}_t^2 - {\bf B}_t^2\right)
+ \frac 12 \sum_a \left( (\partial_t\mbf{\vartheta}_t^a)^2 -
\omega_a^2\mbf{\vartheta}_t^{a2}\right)
+ \sum_a \omega_{pa}\mbf{\vartheta}_t^{a}\cdot {\bf E}_t\right]\ ,
\ee
where the electromagnetic degrees of freedom are described by vector
and scalar potentials, respectively, ${\bf A}_t$ and $\varphi_t$, 
so that ${\bf E}_t = -\mbf{\nabla} \varphi_t -
\partial_t {\bf A}_t$ and ${\bf B}_t = \mbf{\nabla}\times {\bf A}_t$.
The units are chosen in this Section so that $c=1$.
The polarization of the medium is expressed via the matter fields as
${\bf P}_t = \sum_a \omega_{pa} \mbf{\vartheta}_t^a$.
The least action principle for the scalar potential $\varphi_t$ leads 
to the Gauss law, $\mbf{\nabla}\cdot {\bf D}_t = 0$, for the vector
potential ${\bf A}_t$ to the Maxwell's equation, $\partial_t {\bf D}_t 
= \mbf{\nabla}\times {\bf B}_t$, and for the matter field
$\mbf{\vartheta}_t^a$ to the medium polarization evolution equation of the 
Lorentz model with no attenuation, $\gamma_a=0$. 
The second Maxwell's equation and the Gauss law
for the magnetic field follow from the relation 
${\bf B}_t = \mbf{\nabla}\times {\bf A}_t$ by taking its time
derivative and divergence, respectively. The energy of the system
coincides with the canonical Hamiltonian which is obtained by
a Legendre transformation \cite{arnold} 
of the Lagrangian $L$ for the velocities $\partial_t{\bf A}_t$
and $\partial_t \mbf{\vartheta}_t^a$.  Doing the Legendre transformation,
we find the canonical Hamiltonian (energy) of the system
\be
\label{energy}
{E}_t = \frac 12 \int\! d{\bf r} \left[
{\bf E}_t^2 + {\bf B}_t^2 + \sum_a\left(
\mbf{\pi}_t^{a2} + \omega_a^2 \mbf{\vartheta}_t^{a2}\right)\right] = \frac 12
\left|\Psi_t^F\right|^2\ ,
\ee
where $\mbf{\pi}_t^a=\delta L/\delta(\partial_t\mbf{\vartheta}^a)$ 
are canonical momenta of the matter fields. 
To get the last equality in (\ref{energy}), we have used the relations 
$\mbf{\xi}_t^{2a} = \mbf{\pi}_t^a$
and $\mbf{\xi}_t^{2a-1} = \omega_a\mbf{\vartheta}_t^a$ which follow
from comparison of the canonical Hamiltonian equations of motion for
the canonically conjugate variables $\mbf{\vartheta}_t^a$ and $\mbf{\pi}_t^a$
 and Eqs. (\ref{aux1}) and (\ref{aux2}) with $\gamma_a=0$.
Note that the canonical momentum conjugate to the vector potential
${\bf A}_t$ coincides with $-{\bf D}_t = -{\bf E}_t - {\bf P}_t$, not
$-{\bf E}_t$ in this system. Therefore, 
the coupling between the electromagnetic and matter degrees of
freedom is included into the term 
${\bf E}_t^2 = ({\bf D}_t -{\bf P}_t)^2$ of the canonical 
Hamiltonian. Equation (\ref{energy}) becomes the conventional expression for
the electromagnetic energy in a passive medium \cite{landau} when $\mbf{\pi}^a_t$
and $\mbf{\vartheta}^a_t$ are replaced by the corresponding 
solutions of the equations
of motion with initial conditions $\mbf{\pi}^a_0=\mbf{\vartheta}^a_0=0$. 
The energy conservation can be deduced either from the Noether theorem
(because $E_t$ is the Noether integral of motion corresponding to
the time translational symmetry of the action) or directly from the norm
conservation of $\Psi_t^F$ (because the evolution operator $\exp(t{\cal H}^F)$ is
unitary when $\gamma_a=0$).

In numerical simulations, an exact unitary evolution operator ${\cal
U}^Q_{\Delta t}$ is replaced by its approximation ${\cal G}_{\Delta
t}^Q$. However, the energy remains conservative:

{\bf Theorem 6.1}. The split algorithm is unitary for multiresonant
Lorentz models with no attenuation, that is, the split algorithm
preserves the energy $E_{t+\Delta t}=E_t$.

{\bf Proof}. In the field representation, ${\cal H}_0^{F*} = -{\cal H}_0^F$
and ${\cal V}^{F*} = -{\cal V}^F$ and, therefore, ${\cal G}_{\Delta t}^F$ is
unitary. As a result, the algorithm preserves the initial wave packet
energy and the norm,
\be
\label{uni1}
\left|{\cal G}_{\Delta t}^F\Psi_t^F\right| = \left|\Psi^F_{t+\Delta
t}\right|= \left|\Psi_t^F\right|\ .
\ee
In the induction representation, the energy coincides with the norm
of $\Psi_t^I$ in the measure space. The measure is determined by the
transformation law $\Psi_t^I = {\cal S}\Psi_t^F$,
\be
\label{enorm}
{E}_t = \frac 12 (\Psi_t^F,\Psi_t^F) = \frac 12 (\Psi_t^I,\mu
\Psi_t^I)\equiv \frac 12 \left|\Psi_t^I\right|_\mu^2\ ,
\ \ \ \ \ \mu = {\cal S}^{-1*}{\cal S}^{-1}\ .
\ee
Since ${\cal H}^I$ is similar to ${\cal H}^F$,  
the Hamiltonian ${\cal H}^I$ is anti-Hermitian
relative to the $\mu$ scalar product,
\be
\label{mu1}
{\cal H}^{I*}\mu = -\mu{\cal H}^I\ .
\ee
The norm conservation (unitarity) in the split
algorithm requires in addition that
the amplification matrix ${\cal G}_{\Delta t}^I$ satisfies
the unitarity condition
\be
\label{mu2}
{\cal G}_{\Delta t}^{I*}\mu {\cal G}_{\Delta t}^I = \mu\ .
\ee
This is indeed the case. To prove (\ref{mu2}), we show that
${\cal H}_0^I$ and ${\cal V}^I$ satisfy the condition (\ref{mu1})
and, hence, the product of their exponentials is a unitary operator
relative to the $\mu$ scalar product. Consider
${\cal H}^I = {\cal S}{\cal H}^F{\cal S}^{-1} = {\cal H}_0^I 
+ {\cal V}^I$ so that ${\cal H}^F = {\cal S}^{-1}{\cal H}_0^I
{\cal S} + {\cal S}^{-1}{\cal V}^I{\cal S}$.
For the Lorentz model,
\be
\label{mu3}
{\cal S}^{-1}{\cal H}_0^I{\cal S} = 
\pmatrix{{\cal H}_0 & 0\cr 0 & 0} = - 
\left({\cal S}^{-1}{\cal H}_0^I{\cal S} \right)^* \ .
\ee
Therefore ${\cal H}_0^I$ satisfies (\ref{mu1}).
From the anti-Hermiticity of ${\cal H}^F$ and (\ref{mu3})
it follows that 
\be
\label{mu4}
\left({\cal S}^{-1}{\cal V}^I{\cal S}\right)^* = - 
{\cal S}^{-1}{\cal V}^I{\cal S}\ .
\ee
Hence, ${\cal V}^I$ also satisfies (\ref{mu1}). Thus,
\be
\label{uni2}
\left|{\cal G}^I_{\Delta t}\Psi_t^I\right|_\mu =
\left|\Psi_{t+\Delta t}^I\right|_\mu = \left|\Psi_t^I\right|_\mu\ ,
\ee
which completes the proof.

The norm (energy) conservation can be used to control numerical 
convergence, especially when the aliasing problem
in the fast Fourier transform is present, i.e., 
when parameters of the medium are
discontinuous functions in space. In a properly designed algorithm
the loss of energy (norm) due to attenuation should be
controlled by the symmetric part of the Hamiltonian operator
\be
\label{attoper}
\partial_t {E}_t = -\sum_a \gamma_a |\mbf{\xi}_t^{2a}|^2 \equiv
\frac 12\, (\Psi_t^Q, {\cal V}_\gamma^Q\Psi_t^Q) \leq 0\ ,
\ee
where ${\cal V}_\gamma^{Q*} = {\cal V}_\gamma^Q = ({\cal H}^{Q*} 
+ {\cal H}^Q)/2 \leq 0$ (a negative semidefinite operator) 
which is, in this case, a diagonal matrix with nonpositive
elements.

\section{Stability of the algorithm}
\setcounter{equation}0

The norm of an operator ${\cal H}$ is defined as
\be
\label{norm}
\|{\cal H}\| = \sup_{|\Psi|=1} |{\cal H}\Psi|\ .
\ee
If the operator is normal, that is, it commutes with its adjoint, 
then its norm coincides with its spectral radius $\rho({\cal H})$, 
the supremum absolute value
of its eigenvalues. In general, $\rho({\cal H})\leq \|{\cal H}\|$.
A family of amplification operators (matrices) ${\cal G}_{\Delta
t}(\alpha)$ 
is called conditionally stable if 
there exists a constant $C(\tau, T)$ such that \cite{richt}
\be
\label{stability}
\|{\cal G}_{\Delta t}^n(\alpha)\| \leq C(\tau, T)\ ,
\ee
for all $\Delta t \in (0,\tau)$, all $0\leq n\Delta t\leq T$ 
for some positive $\tau$ and $T$, and all parameters $\alpha$.   
The unconditional stability of ${\cal G}_{\Delta t}(\alpha)$ means
that (\ref{stability}) holds uniformly in $n\geq 0$ for any $\Delta t>0$ 
and for all $\alpha$, that is, $C$ is independent of $T$ and $\tau$. 
Parameters $\alpha$ can be all wave vectors
${\bf k}$ supported by the grid or simply grid values of the
position vector ${\bf x}$. They can also include parameters of the
medium. Note that if ${\cal G}_{\Delta t}$ is not normal, then 
$\rho({\cal G}_{\Delta t}) \leq \|{\cal G}_{\Delta t}\|$ and,
therefore, the von Neumann condition $\rho({\cal G}_{\Delta t})\leq 1$
is no longer sufficient for stability, while still being necessary.

{\bf Theorem 7.1}. For multiresonant Lorentz models, the split
algorithm in the field representation is unconditionally stable.

{\bf Proof}. We shall prove that 
\be
\label{fy3}
\|{\cal G}_{\Delta t}^F\| \leq  1\ ,
\ee
which leads to the theorem statement
\be
\label{start}
\|({\cal G}_{\Delta t}^F)^n\| \leq \|{\cal G}_{\Delta t}^F\|^n \leq 1\ .
\ee
By definition and making use of the inequality, $\|{\cal A\, B}\| \leq
\|{\cal A}\|\, \|{\cal B}\|$, we get 
\ba
\nonumber
\| {\cal G}_{\Delta t}^F\| &=& \|e^{\Delta t {\cal H}_0^F/2} e^{\Delta t
{\cal V}^F}
e^{\Delta t {\cal H}_0^F/2}\| \\ \label{block0}
&\leq& \|e^{\Delta t {\cal H}_0^F/2}\|^2
\ea
because $e^{\Delta t {\cal V}^F}$ is a unitary operator, so its norm
equals 1. The operator $e^{t{\cal H}_0^F}$ is block-diagonal 
(see (\ref{freeHF}) and (\ref{eHM2})). The norm
of a block-diagonal operator is the maximal norm of its blocks.
The upper left corner block is given by the unitary operator 
$e^{t{\cal H}_0}$ whose norm equals 1. We have then
\be
\label{block}
\|e^{t{\cal H}_0^F}\|
= \max_a \left\{1,\ \|e^{t {\cal H}_{Ma}^F}\|\right\}.
\ee
The norm of the exponential of ${\cal H}_{Ma}^F$
can be found by direct
calculation using the fact that $\| {\cal A}\|^2 = \|{\cal A}^*
{\cal A}\| = \rho({\cal A}^*{\cal A})$ 
and the explicit form of $e^{t{\cal H}_{Ma}^F}$ given in (\ref{eHM}).
For small attenuation, $\omega_a^2 -\gamma_a^2= \nu_a^2 \geq 0$, 
we define $z_a= (\gamma_a/\nu_a)\sin(\nu_at)$ so that the largest
eigenvalue has the form
\ba
\nonumber
\|e^{t{\cal H}_{Ma}^F}\|^2 &=&\| (e^{t{\cal H}_{Ma}^F})^*e^{t{\cal
H}_{Ma}^F}\| \\
\label{fy}
&= &e^{-2\gamma_at}\left(1 + 2z_a^2 + 2z_a\sqrt{1+z_a^2}\right)\ . 
\ea
Since $z_a \leq \gamma_at \equiv y$ for $t\geq 0$, the function (\ref{fy}) 
is bounded from above by $f(y)=e^{-2y}(1+2y^2+2y\sqrt{1+y^2})$.
It is not hard to verify that the derivative $f^\prime(y)$ is
negative for all $y>0$, and that $f(0)=1$. Hence, replacing $t$ by $\Delta
t/2$, we conclude that
\be
\label{fy2}
\|e^{\Delta t{\cal H}_{Ma}^F/2}\| \leq 1\ , 
\ee
from which (\ref{fy3}) follows.
For large attenuation (like in Drude metals),
$\omega_a^2 - \gamma_a^2 =-\nu^a_a\leq 0$,  in  (\ref{fy}) 
we get
$z_a= (\gamma_a/\nu_a)\sinh(\nu_a t)\equiv z_a(t)$.
For $t\geq 0$ the latter relation defines the inverse
function $t=t(z_a)$. 
Once again, the derivative of (\ref{fy}) with respect to $z_a$
can be shown to be negative for all positive
$z_a$ while at $z_a=0$ the function equals 1. So inequalities
(\ref{fy2}) and (\ref{fy3}) hold in this case too. 
This completes the proof.

The proof of Theorem 7.1 given above 
is not the most economical. However,
the idea of estimating the norm of the exponential of the matter
Hamiltonian in order to investigate stability of the algorithm
can be applied numerically to systems more general than 
the Lorentz model because ${\cal H}_M^F$ is local on the grid, that is,
it does not contain derivatives. So the exponentials of ${\cal H}_M^F$
and its adjoint are not expensive to calculate numerically for some
trial values of $\Delta t$ to see if (\ref{fy2}) holds.

We give an alternative proof of the unconditional stability
in the case of the induction representation of the multi-resonant
Lorentz model where an analytical
expression of the exponent of the matter Hamiltonian
is too hard to find, not to mention its norm. 
We shall make use of the following obvious lemma.

{\bf Lemma 7.2}. Let a vector $\psi_t$, $t\geq 0$,  
be a solution of the linear equation
$\partial_t \psi_t = ({\cal H} +{\cal V})\psi_t$ where the linear operators
${\cal H}$ and ${\cal V}$
satisfy the conditions ${\cal H}^* =-{\cal H}$ and 
${\cal V}^* = {\cal V}\leq 0$ (negative semidefinite). Then $|\psi_t|\leq |\psi_0| $
for all $t\geq 0$.

The proof follows from an obvious relation
\be
\nonumber
\partial_t|\psi_t|^2 = 2(\psi_t,{\cal V}\psi_t)\leq 0\ .
\ee
As a consequence we also get
\be
\label{att}
\|e^{t({\cal H}+{\cal V})}\| \leq 1\ ,
\ee
for all $t\geq 0$.

{\bf Theorem 7.3}. For multiresonant Lorentz models, the split
algorithm in the induction representation is unconditionally stable.

{\bf Proof}. If the attenuation is absent, the amplification matrix
${\cal G}_{\Delta t}^I$ is unitary with respect to the energy scalar
product ($\mu$-scalar product) as is shown in (\ref{uni2}). Hence,
$\|({\cal G}_{\Delta t}^I)^n\|_\mu = 1$. 
When the attenuation is switched on, the 
unitarity of the amplification matrix can get violated only
through $\exp(\Delta t{\cal V}^I)$ because
 the operator $\exp(\Delta t {\cal H}_0^I)$ 
remains unitary with respect to the $\mu$-scalar product. The idea is
to prove the unconditional stability with respect to the $\mu$-norm. 
The theorem statement would follow from
the equivalence of the Euclidean and $\mu$- norms. Recall that
two norms $|\Psi |$ and $|\Psi |^\prime$ are equivalent
if there exist two positive constants $C_{1,2}$ such that
\be
\nonumber
C_1|\Psi | \leq |\Psi |^\prime \leq C_2 |\Psi | \ ,
\ee
for all $\Psi$. All topological properties of the space spanned
by $\Psi$ are the same in one norm as in the other; in particular
convergence of a sequence, boundedness of a set, boundedness of a linear
operator, and uniform boundedness of a family of operators are all
invariant concepts under a change of one norm to the other.
If $|\Psi |^\prime =  |\Psi |_\mu$, then  $\|{\cal A}\|_\mu = 
\|{\cal S}^{-1}{\cal A}{\cal S}\|$ and
\ba
\nonumber
|\Psi |_\mu &=& |{\cal S}^{-1}\Psi| \leq \|{\cal S}^{-1}\|\ |\Psi|\ ,\\
\nonumber
|\Psi | &=& |{\cal S}\Psi |_\mu \leq \|{\cal S}\|_\mu |\Psi|_\mu = 
\|{\cal S}\|\ |\Psi|_\mu\ ,
\ea
so that the two norms are indeed equivalent
\be
\label{normeq}
\|{\cal S}\|^{-1}|\Psi|\, \leq\, |\Psi|_\mu\, \leq\, 
\|{\cal S}^{-1}\|\ |\Psi|\ .
\ee  
Since $\|{\cal A}{\cal U}\| = \|{\cal A}\|$ for any 
bounded operator ${\cal A}$
and a unitary operator ${\cal U}$, we infer that
\be
\label{thind1}
\|({\cal G}_{\Delta t}^I)^n\|_\mu \leq \|{\cal G}_{\Delta t}^I\|_\mu^n
= \|e^{\Delta t {\cal V}^I}\|^n_\mu
= \|e^{\Delta t {\cal S}^{-1} {\cal V}^I{\cal S}}\|^n\ .
\ee
When $\gamma_a=0$ (no attenuation), 
the operator ${\cal S}^{-1}{\cal V}^I{\cal S}$
is skew-symmetric (cf. (\ref{mu4})). When $\gamma_a\neq 0$,
the operator ${\cal S}^{-1}{\cal V}^I{\cal S}$ acquires
an addition which is a diagonal
operator with nonpositive elements as follows from 
(\ref{HI}) and (\ref{lmhm}). Therefore the inequality
(\ref{att}) must
hold for it as a consequence of Lemma 7.2,
$$
\|e^{\Delta t {\cal S}^{-1} {\cal V}^I{\cal S}}\| \leq 1\ ,
$$
from which the uniform boundedness of the family 
$({\cal G}_{\Delta t}^I)^n$ with respect to the $\mu$-norm 
immediately follows.
By the equivalence of the two norms (\ref{normeq}), the 
family $({\cal G}_{\Delta t}^I)^n $ is also uniformly
bounded in the Euclidean norm, 
\be
\label{indbound}
\|({\cal G}_{\Delta t}^I)^n\|
\leq \|{\cal S}\|\,\|{\cal S}^{-1}\|\ ,
\ee
which completes the proof.
 
{\bf Comment}. The same idea of making use of the norm equivalence, 
which actually goes in line with 
the Kreiss matrix theorem (its last part) \cite{kreiss,richt}, 
can be applied to analyze the stability of
the split algorithm for generic passive media. It is not
hard to find a quadratic Lagrangian local in time 
such that the corresponding Euler equations describe
propagation of an electromagnetic pulse in generic
non-absorbing media. Due to time translation symmetry,
the system should have a conserved quantity
according to the Noether theorem \cite{arnold}. This integral of motion 
coincides with the canonical Hamiltonian which is a
quadratic form of $\Psi_t^Q$ if the linear
response approximation is valid.
By analogy with the $\mu$-norm, one could try
to identify the canonical Hamiltonian with 
the new norm of $\Psi_t^Q$ which is conserved by construction
and, hence, in an attenuation-free medium
the corresponding evolution operator is unitary.
Thus, it would always be possible to arrange the split
so that the amplification operator is unitary too.
From the physical point of view, 
it is then naturally expected that, when absorption is
added to the system,  the attenuation operator
${\cal V}_\gamma^Q$ would generally satisfy the condition 
(\ref{attoper}) because Fourier amplitudes of fields 
are exponentially attenuated in passive media. 
The latter would make it possible to apply Lemma 7.2
to prove the unconditional stability of the amplification
operator with respect to the norm defined by the canonical
Hamiltonian along the lines similar to the proof of Theorem 7.3.
An obstacle for this rather natural idea to generalize 
Theorem 7.3 to generic passive media
is that the canonical Hamiltonian is not, in general, positive
definite. It becomes positive only on solutions of the equations of
motion for matter fields, which
is a rather common feature of Lagrangian systems with 
higher order time derivatives. 
Thus, the canonical Hamiltonian does not always define
a positive definite quadratic form in the Hilbert space
for a generic passive media, and, hence, cannot serve 
as a new (conservative) norm. 
The study of conditions on 
attenuation-free media
under which a positive definite and conserved quadratic
form does exist goes beyond the scope of this paper since
it would  require the  canonical formalism and the Noether
theorem for theories with 
higher-order time derivatives, which is rather involved
for generic passive media. The question can be addressed more
easily for each particular medium model of interest. 
However, 
the unconditional stability might be excessive
as far as practical needs are concerned.
It is more important to make the split algorithm
convergent 
for a generic passive medium. Then  one should use
the equivalence of (conditional) stability and
convergence according to the fundamental 
convergence theorem due to Kantorovich \cite{kantor,richt}.    

Our findings in this latter approach are summarized in the following theorem.

{\bf Theorem 7.4}. Suppose that the medium response function 
satisfies the causality conditions (that is, its Fourier transform
has poles only in the lower half of the frequency plane, 
$ {\rm Im}\,\omega\leq 0$). Let ${\cal U}_t$ be an exact evolution matrix
in the grid representation (as defined in Section 3), and 
${\cal G}_{\Delta t}$ be an amplification matrix in some
third order split algorithm. Then for band-width limited
wave packets the split algorithm is (conditionally)
stable and for all $0\leq n \leq N$, $T= N\Delta t$, and $0<\Delta t<\tau$, 
there exist a constant $C_m$, which  
depends only on the medium parameters, and a constant 
$W_m$, which depends also on $\tau$, such that
\ba
\label{th1}
\|{\cal G}_{\Delta t}^n\| &\leq& C_m +\delta(\Delta t, T)\ ,\\
\label{th2}
\delta(\Delta t, T)& =& C_m\left(e^{W_mT\Delta t^2} -1\right)=
O(\Delta t^2)\ ,\\
\label{th3}
\|{\cal U}_{n\Delta t} - {\cal G}_{\Delta t}^n\| &\leq& \delta (\Delta t, T)\ .
\ea

{\bf Remark}. Before proving the theorem, let us discuss its significance
for practical applications. 
Inequality (\ref{th1}) implies conditional stability, while 
(\ref{th3}) establishes a relation between the accuracy (and
convergence) of the 
split approximation and the uniform bound in
the stability condition (\ref{th1}). By making the time step
$\Delta t$ smaller, any desired accuracy can be achieved during
the total (fixed) simulation time $T$. The latter implies, of
course, that the grid is assumed to be chosen fine enough
(in accord with the Shannon sampling theorem) to accurately
reproduce the initial pulse configuration via the fast Fourier
method. Indeed, let $\Psi_{n\Delta t}^{app}= {\cal G}_{\Delta t}^n
\Psi_0$ be a simulated solution, and $\Psi_t= {\cal U}_t\Psi_0$
be an exact solution, then from (\ref{th3}) it follows that
\be
\label{th4}
|\Psi_t -\Psi_t^{app}|\leq \delta(\Delta t, T)|\Psi_0| = O(\Delta
t^2)\ ,
\ee 
for all $0\leq t\leq T$ and any {\it fixed} total simulation time $T$
which is roughly $2L/c$ where $L$ is the simulation box size and
$c$ the speed of light. Now we turn to the proof.

{\bf Proof}. In the grid Fourier basis, ${\cal U}_t\Psi_0^Q = 
\sum_{\bf k} \Psi_t^Q({\bf k}) e^{i{\bf k}\cdot{\bf x}}$,
where ${\bf k}$ spans the dual lattice. By construction of the
Hamiltonian,  each Fourier mode $\Psi_t^Q({\bf k})$ evolves exactly
as in the continuum case. Since the medium response function
satisfies the causality conditions, Fourier amplitudes of the 
electromagnetic and response fields as well as their time
derivatives are bounded functions of time. The amplitudes cannot
grow infinitely large because of dissipation \cite{landau}.
The number of Fourier modes is finite on the grid (only bandwidth
limited initial wave packets are considered) and, hence,
$|\Psi_t^Q| \leq C_Q$ for all $t\geq 0$ because components of the auxiliary
field $\xi_t$ are linear combinations of the response field and
its time derivatives.
The latter inequality 
is equivalent to the evolution matrix being uniformly bounded
for all $t\geq 0$,
\be
\label{bound2}
\|{\cal U}_t\| \leq C_m\ .
\ee
Let ${\cal U}_{\Delta t} - {\cal G}_{\Delta t} = \Delta t^3{\cal
W}_{\Delta t}$ and ${\cal W}_{\Delta t} = {\cal W}_0 + O(\Delta t)$
for small $\Delta t$ according to 
a third order split (cf. (\ref{split}) - (\ref{splitw})). 
Let $W_m = C_m \sup_{\Delta t}\|{\cal W}_{\Delta t}\|$ 
for $0 < \Delta t< \tau$
and some positive finite $\tau$.  Using the semigroup property
${\cal U}_{\Delta t}^k = {\cal U}_{k\Delta t}$ and (\ref{bound2}) 
we infer that
\ba
\nonumber
\|{\cal G}_{\Delta t}^n\| &=&\left\|{\cal U}_{\Delta t}^n - 
\left({\cal U}_{\Delta t}^n - {\cal G}_{\Delta t}^n\right)\right\|\\
\label{bound3b}
&\leq & \|{\cal U}_{\Delta t}^n\| + 
\|{\cal U}_{\Delta t}^n -{\cal G}_{\Delta t}^n\|\\
\label{bound3a}
&\leq & C_m + \|{\cal U}_{\Delta t}^n -({\cal U}_{\Delta t} - 
\Delta t^3 {\cal W}_{\Delta t})^n\|\\
\label{bound3}
&=& C_m + \| -\Delta t^3 \left(
\sum_{k=0}^{n-1}{\cal U}_{\Delta t(n-k-1)}{\cal W}_{\Delta t}
{\cal U}_{k\Delta t}\right) + \cdots\|\\
\label{bound4}
&\leq& C_m + C_m\left[\left(1+\Delta t^3 W_m)^n - 1\right)
\right]\\
\label{bound5}
&\leq& C_m + \delta(\Delta t, T)\ .
\ea
Inequality (\ref{th3}) readily follows
from comparing the right hand side of (\ref{bound3b}) with 
those of (\ref{bound3a})-(\ref{bound5}). This completes the proof.

\section{Convergence and accuracy analysis}
\setcounter{equation}0
To estimate the accuracy of the algorithm at a fixed finite grid size $N$, 
consider the following
quantity
\be
\label{accb1}
\beta_n(N,\Delta t)= \left|\left(({\cal U}_{\Delta t}^Q)^n - ({\cal G}_{\Delta
t}^Q)^n\right)
\Psi_0^Q\right|/\left|\Psi_0\right| \leq 
\left\|({\cal U}_
{\Delta t}^Q)^n - ({\cal G}_{\Delta }^Q)^n\right\|
\ee
which specifies a deviation of the approximate solution
from the exact one relative to a given norm. Here ${\cal U}_{\Delta t}^Q$
is an exact evolution operator.
The accuracy estimate $\beta_n(N,\Delta t)$
is a norm dependent quantity. The choice of norm is usually
determined by practical needs. We use the norm related to 
the electromagnetic energy of the system and investigate, first,
the behavior of $\beta_n(N,\Delta t)$ as $\Delta t$ goes to zero,
while $\Delta t n =t $ remains fixed and does not exceed some positive
constant, $t\leq T$. 

{\bf Theorem 8.1}. For multi-resonant Lorentz models, there 
exists a positive constant $W^Q$ such that
\be
\label{th41}  
\left\|({\cal U}_
{\Delta t}^Q)^n - ({\cal G}_{\Delta }^Q)^n\right\|
\leq \Delta t^2 TW^Qc_Q^2\ ,
\ee
where $c_F=1$ and $c_I = \|{\cal S}\|\, \|{\cal S}^{-1}\|$, for 
all $\Delta t \in (0,\tau)$ and $n\Delta t \leq T$.

{\bf Proof}. In the field representation $Q=F$, 
$\|({\cal U}_{\Delta t}^F)^n\| \leq 1$ and
$\|{\cal G}_{\Delta t}^F)^n\|\leq 1$ for any integer $n$, as a 
consequence of Lemma 7.2 for the multi-resonant Lorentz model.
The same inequalities hold in the induction representation
if the norm is replaced by the $\mu$-norm.
According to the split algorithm (\ref{split})-(\ref{splitw}), 
${\cal U}_{\Delta t}^Q -
{\cal G}_{\Delta t}^Q = \Delta t^3 {\cal W}_{\Delta t}^Q$. 
Let $W^Q = 
\sup_{\Delta t} \|{\cal W}_{\Delta t}^Q\|$ for $\Delta t
\in (0,\tau)$ for some positive $\tau$ (a maximal
time step used in simulations). We then have the following
chain of inequalities that lead to the theorem statement
\ba
\label{accb2}
\left\|({\cal U}_
{\Delta t}^F)^n - ({\cal G}_{\Delta }^F)^n\right\| &=&
\left\|\sum_{k=1}^{n-1}({\cal U}_{\Delta t}^{F})^k\left(
{\cal U}_{\Delta t}^F -{\cal G}_{\Delta t}^F\right)
({\cal G}_{\Delta t}^F)^{n-k}\right\| \\
\label{accb3} 
&\leq& (n-1)\Delta t^3 \|{\cal W}_{\Delta t}^F\| \\
\label{accb3a}
&\leq& \Delta t^2 TW^F\ ,
\ea
In the case of the induction representation, inequality (\ref{accb3})
holds relative to the $\mu$-norm. The theorem statement (\ref{th41}) follows
from the norm equivalence (\ref{normeq}), 
$c_I^{-1}\|{\cal A}\|\leq \|{\cal A}\|_\mu \leq c_I \|{\cal A}\|$
for any operator ${\cal A}$. The proof is complete.

{\bf Remark}. In simulations, the continuum limit $N\rightarrow
\infty$ is never achieved.
Hence the operators in the split algorithm (\ref{kato}) remain bounded
versus the unbounded case of (\ref{kato}). 
It is known that the convergence rate of $\beta_n(\infty,
\Delta t)$ as $\Delta t\rightarrow 0$ estimated by the operator 
norm as in the right hand side of (\ref{accb1}) is no longer of
order $O(\Delta t^2)$ but rather of $O(\sqrt{\Delta t})$ 
(see, e.g., \cite{ichi} and references therein). For
unbounded operators, the estimate (\ref{accb3a}) is not valid.
This suggests that the convergence rate $\beta_n(N, \Delta t)$
may depend, even significantly, on the initial vector $\Psi_0$
as $N$ increases.  

In a general case, the quantity $\delta (\Delta t, T)$ 
in Theorem 7.4 determines the accuracy of the split algorithm
with respect to the norm (\ref{norm}) on a finite grid. 
To make simulation errors
small, it is sufficient to require that
\be
\label{acc1}
\|{\cal U}_{\Delta t} - {\cal G}_{\Delta t}^Q\|=
\Delta t^3\| {\cal W}_{\Delta t}^Q\| <\!\! < 1\ .
\ee
Making use of (\ref{splitw}) and the fact that the norm of a matrix
does not exceed the maximal norm of its blocks, we infer
for a multi-resonant Lorentz model
that, in order for (\ref{acc1}) to hold, the following inequalities
are sufficient:
\be
\label{acc2}
\omega_{pa}\Delta t <\!\! < 1\ ,\ \ \ \ 
\omega_{max}\Delta t <\!\! < 1 \ ,\ \ \ \omega_a\Delta t <\!\! < 1\ ,\ \ \ \ 
\gamma_a\Delta t <\!\! < 1\ ,
\ee
and, yet another one, 
\be
\label{acc3}
\frac{|\mbf{\nabla }\omega_{pa}|}{\omega_{pa}}\ c\Delta t <\!\! < 1\ .
\ee
Here $\omega_{max} $ is the maximal frequency of the initial wave packet. 
The right hand side of (\ref{acc1}) is a sum of two types of
terms. There are terms  
of the cubic order in numbers (\ref{acc2}) as well as a term linear in (\ref{acc3}) with
the coefficient being quadratic in (\ref{acc2}). 
The ratio in (\ref{acc3}) can roughly be estimated from
$|\mbf{\nabla} \omega_{pa}| \leq \omega_{pa}/\Delta { r}$  with $\Delta {r} $ 
being  the grid step.
The condition (\ref{acc3}) implies then that the distance 
traveled by the wave packet during one time step should be
much smaller than the grid step. 

To complete the discussion, one should also analyze the accuracy
of the Gauss law (\ref{gauss}).
Note that the constraints are automatically fulfilled
in the continuum theory due to the Dirac involution relations (\ref{dirac}).
By projecting the continuum theory onto a finite grid and replacing
the exact evolution operator by its approximation in the split algorithm, the
involution relations might be violated, thus leading to errors and
potential instabilities of the algorithm. A good example
of this kind is numerical general relativity (although the
nonlinearity of the equations of motion plays the central
role in generating instabilities due to the violation of the
Dirac involution relations).

It is not hard to be convinced that the Gauss law (\ref{gauss}) is
equivalent to the following constraint on state vectors 
\be
\label{gausslaw1}
{\cal C}^{I} \Psi_t^I = 0 \ ,\ \ \ \ 
{\cal C}^{I} =\pmatrix {{\cal C} & 0\cr 0 & 0}\ ,\ \ \ {\cal C} = 
\pmatrix { 1 & 0\cr 0 & 1} {\cal P}_\|\ ,\ \ \ {\cal  P}_\| = 1 -{\cal P}_\bot\ .
\ee
The operator ${\cal P}_\|$ projects a vector field onto its
longitudinal component. In other words, it acts as the identity
operator if the vector field is conservative, and it annihilates
any rotational vector field (which is the curl of another vector field).
In the field representation we get ${\cal C}^F = {\cal S}^{-1}{\cal C}^I
{\cal S}$ with ${\cal S}$ defined in (\ref{ftoi}).
On the grid, the action of the operator ${\cal C}^Q$ is defined by
the rule (\ref{grid5}), that is, by (\ref{gausslaw1}) we understand
$({\cal F}{\cal C}^Q{\cal F}^*) {\cal F}\Psi_t^Q=0$.
Thus, the Gauss law
requires that the Fourier transform of the  inductions should not
acquire components parallel to wave vectors of the dual grid.
This is obviously guaranteed if the exact evolution operator, 
${\cal U}_t^Q = \exp(t{\cal H}^Q)$,
is used to generate the time evolution because
\be
\label{gausslaw2}
{\cal C}^Q {\cal H}^Q =0\ ,
\ee
and, hence, ${\cal C}^Q{\cal U}_t^Q\Psi_0^Q = {\cal C}^Q\Psi_0^Q=0$.
A problem may arise when the approximate evolution
operator, $({\cal G}_{\Delta t}^Q)^n$, is used to evolve
the initial wave packet $\Psi_0^Q$. From linearity of the system,
it is natural to expect that the Gauss law violation
should be of the same order as the accuracy of a numerical solution
of dynamical Maxwell's equations. However, we shall take
a closer look at the problem and find a pleasant result
important in practice, which is stated in the following
theorem.

{\bf Theorem 8.2}. Assuming linear response theory for any passive
medium, the Gauss law 
holds exactly in the split algorithm in the induction
representation.

{\bf Proof}. In the induction representation, identity
(\ref{gausslaw2}) is equivalent to two identities
for the blocks of ${\cal C}^Q{\cal H}^Q$, namely,
${\cal C}{\cal H}_0=0$ and ${\cal C} {\cal V}_{IM} = 0$.
The first one is obvious. The second one follows from
(\ref{VIM}) established for any passive medium.
The key observation is that the identity
\be
\label{gausslaw3}
{\cal C}^I{\cal H}_0^I=0
\ee
holds thanks to the two above identities and 
 (\ref{isplit}). Indeed, in the Fourier basis
(\ref{gausslaw3}) is equivalent to the vanishing of
the triple vector product ${\bf k}\cdot({\bf k}
\times {\bf A})$ for some ${\bf A}$ regular
at ${\bf k}=0$.
Then
from (\ref{gausslaw1}) and (\ref{gausslaw3})
it follows that
\be
\label{gausslaw4}
{\cal C}^I{\cal V}^I = {\cal C}^I({\cal H}^I
-{\cal H}_0^I) =0\ .
\ee
As a consequence of (\ref{gausslaw3}) and (\ref{gausslaw4}),
we infer that
\be
\label{gausslaw5}
{\cal C}^I\left({\cal G}_{\Delta t}^I\right)^n\Psi_0^I = 
{\cal C}^Ie^{\Delta t{\cal H}_0^I/2}e^{\Delta t{\cal V}^I}
e^{\Delta t{\cal H}_0^I/2}
\left({\cal G}_{\Delta t}^I\right)^{n-1}\Psi_0^I
={\cal C}^I\left({\cal G}_{\Delta t}^I\right)^{n-1}\Psi_0^I
={\cal C}^I\Psi_0 =0\ ,
\ee
which is the statement of the theorem.

In the field representation the Gauss law can be enforced 
by means of the projection formalism discussed in Section 1.
The projection operator is, obviously, ${\cal P}=1-{\cal C}^F$.
Its action is computed in the grid representation by the 
fast Fourier method according to the rule (\ref{grid5}).
Without the use of the projection formalism,
the accuracy of the Gauss law
is stated in the following technical proposition.

{\bf Proposition 8.3}. Let $W = \|{\cal W}_{\Delta t}^F\|$
and $W_C = \|[{\cal C}^F, {\cal W}_{\Delta t}^F]\|
C_m(1 + \delta (\Delta t, T))$ where $C_m$ and 
$\delta(\Delta t ,T)$ are defined in Theorem 7.4, then
\be
\label{gausslaw6}
\left|{\cal C}^F\left({\cal G}_{\Delta t}^F\right)^n\Psi_0^F\right|
/|\Psi_0^F| \leq TW_C \Delta t^2 +\Delta t^4 WW_Ce^{TW\Delta t^2}(\Delta
t^2 + T^2/2) = O(\Delta t^2)\ ,
\ee
for all $0\leq n \leq N$, $T=N\Delta t$ and any positive $\Delta t$.

{\bf Proof}. Since ${\cal C}^F{\cal U}_t^F = {\cal C}^F$, assuming
that the initial state $\Psi_0^F$ satisfies the Gauss law we get
\ba
\nonumber
{\cal C}^F({\cal G}_{\Delta t}^F)^n \Psi_0^F &=& 
-{\cal C}^F\left\{({\cal U}_{\Delta t}^F)^n - 
({\cal G}_{\Delta t}^F)^n\right\} \Psi_0^F\\
\nonumber
&=&-{\cal C}^F\sum_{k=1}^{n-1}({\cal U}_{\Delta t}^F)^k
\left({\cal U}_{\Delta t}^F - {\cal G}_{\Delta t}^F\right)
({\cal G}_{\Delta t}^F)^{n-k}\Psi_0^F\\
\label{prop1}
&=& -\Delta t^3 [{\cal C}^F,{\cal W}_{\Delta t}^F]
\sum_{k=1}^{n-1}({\cal G}_{\Delta t}^F)^{n-k}\Psi_0^F 
+\Delta t^3 {\cal W}_{\Delta t}^F \sum_{k=1}^{n-1}
{\cal C}^F({\cal G}_{\Delta t}^F)^{n-k}\Psi_0^F\ .
\ea
Denoting the left hand side of (\ref{gausslaw6}) by $\alpha_n$, we
infer from (\ref{prop1}), by taking the norm of both sides, that 
$$
\alpha_n \leq \Delta t^3 \left\|[{\cal C}^F,{\cal W}_{\Delta
t}]\right\|
\sum_{k=1}^{n-1} \|({\cal G}_{\Delta t}^F)^{n-k}\| +
\Delta t^3 \|{\cal W}_{\Delta t}^F\|\sum_{k=1}^{n-1}\alpha_{n-k}\ ,
$$
for $n>1$ and $\alpha_1 \leq W_C\Delta t^3$.
By Theorem 7.4, powers of the amplification matrix 
${\cal G}_{\Delta t}^F$ are bounded. Hence the following
recursion inequality holds
\be
\label{prop2}
\alpha_n \leq (n-1)\Delta t^3 W_C + \Delta t^3 W (\alpha_{n-1}
+\alpha_{n-2} +\cdots + \alpha_1)\ .
\ee
Iterating (\ref{prop2}) $n-1$ times, we deduce that
\ba
\nonumber
\alpha_n &\leq& (n-1)\Delta t^3 W_C +
\Delta t^3W\left\{(n-2)\Delta t^3W_C +(1+\Delta t^3W)
(\alpha_{n-2}+\alpha_{n-3} +\cdots \alpha_1)\right\}\\
\nonumber   
&\leq & (n-1)\Delta t^3W_C + \Delta t^6WW_C\left\{
\sum_{k=0}^{n-2}(n-2-k)(1+\Delta t^3W)^k + (1+\Delta t^3W)^{n-2}
\right\}\ .
\ea
One can find an explicit form for  the sum in the latter equation. However,
it is a cumbersome expression. For practical purposes,
we give a simpler estimate which is stated in (\ref{gausslaw6}).
First, factor out $(1+\Delta t^3W)^{n-2}$ in the brackets,
and then use obvious inequalities $(1+\Delta t^3W)^{-k}\leq 1$
and $ (1+\Delta t^3W)^{n-2}\leq \exp(TW\Delta t^2)$, which leads to
(\ref{gausslaw6}). 

In the case of the Lorentz model, $W_C=\|[{\cal C}^F,{\cal W}_{\Delta
t}^F]\|$
because all powers of the amplification matrix are uniformly bounded
by 1. For small $\Delta t$, a good estimate can be obtained by 
computing $W_C$ for $\Delta t=0$ using (\ref{splitw}).

The convergence rate as the number of grid points $N$
increases is determined by the convergence rate of 
the fast Fourier transform which is exponential
versus polynomial in finite difference schemes, provided
parameters of the medium are smooth functions of position \cite{psm,fft}.
As is well known from Fourier analysis, the convergence rate 
can be affected for functions which 
have discontinuities \cite{fft}. The latter is, unfortunately, the case
in electromagnetic scattering problems. Suppose there is an interface
between two media. It can be deduced from 
the Maxwell's equations that the components of 
the electric and magnetic fields, ${\bf E}_t$ and ${\bf H}_t$,
tangential to the interface must be continuous, provided there is
no surface electric current on the interface. From the Gauss law
it follows that the components of the inductions, ${\bf D}_t$
and ${\bf B}_t$, normal to the interface must be continuous, provided
there is no surface charge on the interface. In contrast, 
the normal components
of the fields and the tangential components of the inductions can be
discontinuous. Their discontinuities are proportional to 
discontinuities of medium parameters (e.g., discontinuities in
plasma frequencies in Lorentz models).
Therefore, in either the
induction or field representation, there are components which suffer
discontinuities at the interface. Consequently, the convergence rate
of the split algorithm for Maxwell's theory might be 
slower than that in quantum mechanics with
a discontinuous potential because in the latter case 
the wave function remains continuous.

Another source of errors that affects the convergence rate as $N$ increases
is the aliasing problem in the fast Fourier transform on the grid.
Note that, even though the initial wave packet is band-width limited
and the grid is chosen fine enough to eliminate errors in doing its
fast Fourier transform back and forth, the wave packet looses 
this property after the operator $\exp(\Delta t{\cal V}^Q) $ is
applied to it. As a result, the aliasing problem arises in spatial domains
where ${\cal V}^Q$ varies (typically at interfaces between different
types of media). 

The above two problems that also reduce the accuracy of the
algorithm are well known and studied 
in the theory of the fast Fourier transform \cite{fft}. 
The only way to cope with them is to make the grid finer in the areas
where medium parameters have discontinuities. However, the fast Fourier
algorithm requires a uniform equispaced lattice, which might lead
to wasting computer resources if the increased resolution is necessary
only in relatively small areas of the computational volume of the
problem (e.g., only near an interface between two media). There are 
several ways to modify the algorithm when the above problems 
are too expensive to overcome by making a uniform grid finer.

First, the grid can be made fine enough so that the action of powers of the 
Hamiltonian ${\cal H}^Q$ on the state vector $\Psi_0^Q$ is
sufficiently accurate in the Fourier basis as specified by the 
rules (\ref{grid4}) and (\ref{grid5}). The operator
${\cal H}^Q$ is projected onto 
the Krylov space spanned by vectors $({\cal H}^Q)^k\Psi_0^Q$,
$k=0,1,...,n$. Its exponent (the evolution operator) 
is then computed by diagonalizing ${\cal H}^Q$ instead
of using the Lie-Trotter formula. Usually, it is sufficient 
to take a low dimensional Krylov space. This method is known
as the Lanczos method \cite{lanczos}. A detailed study of the Krylov-Lanczos
method as well as other similar pseudospectral methods in
Maxwell theory will be given elsewhere. 

Second, one can give up a uniform grid, while preserving basic
advantages of pseudospectral
methods such as, e.g., exponential convergence.
A possible way to emulate a non-uniform grid in a multiscale
problem is to use wavelet bases. The problem here is to 
compute the action of $\exp(\Delta t{\cal H}_0^Q)$ in the split
algorithm  because
the derivative operator $\mbf{\nabla}$ is not diagonal in
this basis (in contrast to the Fourier basis). However,
${\cal H}_0^Q$ is expected to be sparse in a wavelet basis
so that its direct diagonalization might not be expensive, and
a significant reduction of computational costs can be achieved
in the split algorithm, by using the fast wavelet transform, 
as compared to that in the Fourier basis.
Otherwise, the use of (\ref{inv}) might be helpful in place
of the direct diagonalization method. This approach has proved
to be successful in solving multiscale initial value problems
for the Schr\"odinger equation \cite{daubl}. In the framework of Maxwell's
theory for passive media, additional studies of several issues
in time domain wavelet based algorithms,
like, e.g., stability, would still be needed.

Third, the fast Fourier transform algorithm remains in place
but is applied to an auxiliary uniform grid that is related to 
a non-uniform grid in physical coordinates by a change of variables.
Consider a change of variables ${\bf y}={\bf y}({\bf x})$.
A uniform grid in the new variables ${\bf y}$ would generate
a non-uniform grid in the original Euclidean (physical) coordinates ${\bf x}$.
A desired  local density of grid points in the physical space, to enhance the sampling
efficiency in designated regions, can be achieved by an appropriate 
choice of the functions ${\bf y}({\bf x})$ \cite{cc}. By necessity, the auxiliary
grid spans a rectangular box (with periodic boundary conditions). 
Its pre-image in the physical space would not be a box in general,
save for the case when the map ${\bf y}({\bf x})$ splits into
three individual one-dimensional maps $y_j=y_j(x_j)$, $j=1,2,3$.  
Since, the derivatives are transformed as $\mbf{\nabla}_{\bf x}=
{\cal A}({\bf y})\mbf{\nabla}_{\bf y}$ where the $3\times 3$ matrix
${\cal A}$ is position dependent, the operator ${\cal H}_0^Q$ cannot
be kept in the exponential. The action of its exponential on the state 
vector can be approximated by the leapfrog method in which only the action
of ${\cal H}_0^Q$ on $\Psi^Q$ is required. The latter can be
done by the fast Fourier method according to the rules (\ref{grid5})
and (\ref{grid4}) applied to an operator being a product of position
and derivative dependent operators. 
In contrast to the well studied 
quantum mechanical case, the algorithm appears to be unstable for
media with absorption. In Section 9 a modification of the leapfrog
scheme is proposed to achieve conditional stability.

\section{The temporal leapfrog scheme}
\setcounter{equation}0

Here we discuss a temporal finite difference scheme applied to the Maxwell
theory for passive media in the Hamiltonian formalism. 
As has been pointed out, such a scheme
might be helpful for reducing computational costs by using non-uniform
grids in combination with some pseudospectral methods (e.g., wavelet
bases or the fast Fourier method with a change of variables).
A temporal finite difference scheme can be obtained by the following
procedure. Let us integrate (\ref{sch})
 over the interval $(t, t+n\Delta t)$. We have 
\be
\label{fds}
\Psi_{t+n\Delta t} = \Psi_t + {\cal H} \int\limits_t^{t+n\Delta t} d\tau \Psi_\tau\ 
= \Psi_t + \Delta t {\cal  H} 
\left(\sum_{k=0}^{n-1} C_k^{(n)} \Psi_{t+k\Delta t}\right) + O(\Delta t^{n+1})\ ,
\ee
where the  coefficients $C_k^{(n)}$ used to approximate the integral are well  known for any $n$
as well as the accuracy of the approximation. 
For example, one can use the $3/8$ Simpson rule for $n=3$ or Bode's rule for $n=4$.
The iterating
scheme allows one to compute the wave function at the sequential moment of time if it is known for $n$
preceding moments of time. Only the simplest case $n=2$, for which
the mid-point approximation for the integral 
is taken, leading to $C_0^{(2)}= 0$ and $C_1^{(2)}= 2$, will be
considered in detail.  
It is also known as the leapfrog scheme:
\be
\label{lf1}
\Psi_{t+\Delta t} = \Psi_{t-\Delta t} + 2\Delta t {\cal H} \Psi_t\ .
\ee
The action of the Hamiltonian is computed in a suitable basis
(as has been noted above). 
Apart from violation of the dispersion 
relation of electromagnetic waves, temporal finite difference schemes 
would generally be unstable in media with absorption, in contrast to
the quantum mechanical case. The reason is that 
the Hamiltonian in (\ref{lf1}) is
not anti-Hermitian. Consequently, convergence to the 
continuum solution  would also be violated.
 
A general solution to (\ref{lf1}) can be written in the form
\ba
\label{l1}
\Psi_{n\Delta t} &=& \left({\cal G}_{\Delta t}^{(+)}\right)^n \Psi_{+} +
 \left({\cal G}_{\Delta t}^{(-)}\right)^n \Psi_{-}\ ,\\
\label{l2}
{\cal G}_{\Delta t}^{(\pm)} &=& {\cal H}\Delta t \pm 
\sqrt{1 +{\cal H}^2\Delta t^2} \ ,
\ea
for some initial state vectors $\Psi_0$ and $\Psi_{\Delta t}$ (
the vectors $\Psi_\pm$ are determined by them). Stability
requires that there exists a positive constant $C$ such that
\be
\label{l3}
\left|\Psi_{n\Delta t}\right| \leq C \left(
\left|\Psi_+\right| + \left|\Psi_-\right|\right)\ ,
\ee
for all $0\leq n\leq N$, $T=N\Delta t$ and $0<\Delta t<\tau$.
Note that in general a solution of (\ref{sch}) may have a
legitimate exponential growth if the hermitian part of
the Hamiltonian, ${\cal H}+{\cal H}^*$, is not negative
semidefinite
(see Lemma 7.2). For this reason, a typical stability criterion
would be equivalent to the condition \cite{richt} that there exists some
positive constant $K_1$ such that 
$\|{\cal G}^{(\pm)}_{\Delta t}\| \leq 1 + K_1\Delta t$ uniformly
for all parameters of ${\cal G}^\pm_{\Delta t}$ and for
$0<\Delta t<\tau$, which is clearly the case for
(\ref{l2}) if  ${\cal H}$ is bounded. The latter leads to 
\be
\label{l4}
\left\|\left({\cal G}_{\Delta t}^{(\pm)}\right)^n\right\| \leq 
e^{K_1T}
\ee
so that  a legitimate exponential growth of the solution is allowed,
i.e., $C\sim e^{K_1T}$ in (\ref{l3}).
For passive media, the Hamiltonian satisfies the conditions
of Lemma 7.2 and, hence, no legitimate exponential growth should be
present in a numerical solution in order to achieve convergence.
However, as we shall see shortly, 
the scheme (\ref{lf1}) always generates an exponentially
growing solution for media with attenuation.

Let complex numbers $z=Re^{i\varphi}$ 
be eigenvalues of ${\cal H}\Delta t$.
Since the spectral radius of ${\cal G}_{\Delta t}^{(\pm)}$
does not exceed its norm, the necessary (von Neumann) condition
to suppress an exponential growth of the solution (\ref{l1})
reads
\be
\label{l5}
\rho({\cal G}_{\Delta t}^{(\pm)}) = \max_{z\in D}\left|
z\pm \sqrt{1+z^2}\right| \leq 1\ .
\ee
The aim is to analyze the domain $D$ of the complex plane for which (\ref{l5}) 
holds. Let $\eta = \sqrt{1+z^2}$ and $|\eta|=r$. The two 
inequalities in (\ref{l5}) require that for $z\in D$,
$R^2 + r^2 \pm \xi \leq 1$, where $\xi = \bar{z}\eta +z\bar{\eta}$.
By combining the latter inequalities, one gets $R^2 +r^2\leq 1$
or $r^4 \leq (1-R^2)^2$. On the other hand, 
$r^4 = 1 +R^4 + 2R^2\cos(2\varphi)$. Hence,
$\cos(2\varphi) \leq -1$ which is only possible
if $\varphi = \pm \pi/2$ or $z=\pm iR$. 
The necessary condition (\ref{l5}) is satisfied if 
\be
\label{vn}
\varphi = \pm \pi/2\ ,\ \ \ \ \ 
R^2\leq 1\ .
\ee
This does not yet guarantee that there is no norm growth.
A norm growth, which is polynomial in time, can still occur.

Let us investigate general properties of the solution of
(\ref{sch}) when the Hamiltonian satisfies the von Neumann
condition (\ref{vn}). For any matrix ${\cal H}$ there
exists a similarity transformation so that
${\cal S}^{-1}{\cal H}{\cal S}$ has the Jordan normal
form. Let $h_z$ be a block of the Jordan normal form
corresponding to an eigenvalue $z$ of ${\cal H}$.
Any block $h_k$ is a $q_z\times q_z$ bi-diagonal matrix, 
$q_z\geq 1$, with all the elements
of the diagonal being equal to $z$ and all the elements
on the upper superdiagonal being equal to one.
For $q_z=1$, $h_z =z$ is just a complex number.
The norm of any solution of (\ref{sch}) 
cannot grow faster than $\|\exp(t{\cal H})\|$.
Let a $q_z$-dimensional vector $\phi_t$ satisfy the equation
$\partial_t \phi_t = h_z\phi_t$. For a generic
initial condition, the solution norm grows polynomially, 
$|\phi_t| = O(t^{q_z-1})$ as $t\rightarrow \infty$.
Using the similarity transformation ${\cal S}$, we 
define the corresponding $\mu$-norm of state vectors
and the corresponding matrix norm (cf. (\ref{enorm}))
by setting $\mu = {\cal S}^{-1*}{\cal S}^{-1}$.
From the equivalence of the norms
$\|\cdot\|$ and $\|\cdot\|_\mu$ 
(see the proof of Theorem 7.3), the norm growth
cannot be faster than
\be
\label{jform}
\|e^{t{\cal H}}\|_\mu = \|e^{t{\cal S}^{-1}{\cal H}{\cal S}}\|=
\max_z \|e^{th_z}\| = O(t^{q-1})\ ,\ 
\ \ \  \ q=\max_z q_z\ , \ \ \ t\rightarrow \infty\ ,
\ee
provided $z=\pm iR$. However, a state vector norm
growing polynomially with time  
is unacceptable from the physical
point of view because in any passive medium there is
no physical mechanism for such amplification of the field
amplitudes in the large time limit. Consequently,
we demand that any model Hamiltonian for a
passive medium should be similar to a diagonal 
matrix (i.e., ${\cal H}$ is diagonalizable).
In this latter case, the blocks $h_z$ of the Jordan
normal form of ${\cal H}$ are just complex numbers $z$.
Hence $\|\exp(th_z)\|=|\exp(\pm itR)|=1$ so that the $\mu$-norm
of any solution of (\ref{sch}) is conserved according to
(\ref{jform}). 

Two important conclusions about the leapfrog scheme (\ref{lf1}) 
follow from our analysis.
First, the von Neumann condition (\ref{vn}) is also sufficient
for stability. Indeed, if (\ref{vn}) holds then
$\|{\cal G}_{\Delta t}^{(\pm)}\|_\mu
= \|{\cal S}{\cal G}_{\Delta t}^{(\pm)}{\cal S}^{-1} \|
=\rho({\cal G}_{\Delta t}^{(\pm)}) =1$ 
and, hence, $\|({\cal G}_{\Delta t}^{(\pm)})^n\|_\mu\leq 1$ uniformly
in $n\geq 0$. By the norm equivalence, 
$\|({\cal G}_{\Delta t}^{(\pm)})^n\|$
is also bounded uniformly in $n\geq 0$. 
Second, reversing the argument, we conclude from the 
norm conservation of the stable leapfrog solution that
no attenuation can be added to the Hamiltonian 
without destroying the stability and, consequently,
the convergence to the continuum solution.
Whenever the attenuation is added, the leapfrog solution
would always contain an exponentially growing component,
while this would not be so for a continuum solution
(see Lemma 7.2).
  
Since ${\cal G}_{\Delta t}^{(+)}{\cal G}_{\Delta t}^{(-)}=1$,
only one of the two
independent solutions in (\ref{l1}) would grow exponentially
whenever the attenuation is added.
Theoretically, for ${\cal H}+{\cal H}^*\leq 0$ the exponentially
growing solution can be eliminated by choosing the
initial condition so that
$\Psi_-=0$ which is equivalent to the initial condition 
$\Psi_{\Delta t} = {\cal G}_{\Delta t}^{(+)}\Psi_0$.
Practically, this is never possible due to rounding errors
and/or numerical errors in computing ${\cal G}_{\Delta t}^{(+)}
\Psi_0$.
Even for a  small $|\Psi_-|$ in (\ref{l1}) 
the growing part would eventually
become comparable with the exponentially attenuating 
solution generated by $\Psi_+$. A reduction of the time
step would not be helpful since the constant $K_1$
in (\ref{l4}) is independent of $\Delta t$ while the simulation
time $T$ is fixed by the dimension of the simulation volume
and the speed of light. One needs at least to modify the 
scheme so that there exists a constant $K_p$ such that 
\be
\label{l6}
\|{\cal G}_{\Delta t}^{(\pm)}\| \leq 1 + K_p \Delta t^p\ , \ \ \ p>1\ ,
\ee
for $0<\Delta t<\tau$.
Indeed, it follows from (\ref{l6}) that 
$\|({\cal G}_{\Delta t}^{(\pm)})^n
\|\leq \exp(K_pT\Delta t^{p-1})= 
1 + O(\Delta t^{p-1})$ for all $0\leq n\leq N$ where
$N\Delta t =T$. The norm growth could be reduced 
as much as desired by making
the time step smaller. Next we show how to modify the leapfrog
scheme to make (\ref{l6}) valid for at least $p=3$ and, if the
Hamiltonian is normal, an even stronger result holds, namely, 
$K_p=0$. 
 
 Let ${\cal H}=
{\cal H}_0 + {\cal V}$ where ${\cal V}^*+{\cal V}\leq 0$
(negative semidefinite)
and ${\cal H}_0^* = -{\cal H}_0$. In (\ref{sch}) we make
a substitution $\Psi_t = \exp(t{\cal V})\Phi_t$. The 
new state vector $\Phi_t$ satisfies an equation
with a time dependent Hamiltonian,
\be
\label{scht}
\partial_t \Phi_t = e^{-t{\cal V}}{\cal H}_0e^{t{\cal V}}\Phi_t \equiv
{\cal H}_t\Phi_t\ ,
\ee
and with the same initial condition $\Phi_0=\Psi_0$. Applying 
the leapfrog method to (\ref{scht})   we get 
$\Phi_{t+\Delta t} = \Phi_{t-\Delta t} + 2\Delta t{\cal H}_t\Phi_t$
valid up to $O(\Delta t^3)$. Returning to the initial variables,
we arrive at the following recurrence relation
\be
\label{lfnew}
\Psi_{t +\Delta t} = {\cal L}_{2\Delta t}\Psi_{t-\Delta t} +
2\Delta t {\cal L}_{\Delta t} {\cal H}_0 \Psi_t\ ,
\ee
where ${\cal L}_{\Delta t} = \exp(\Delta t {\cal V})$.  
All the derivative operators are included into the anti-Hermitian
part ${\cal H}_0$ of the Hamiltonian ${\cal H}$, 
while the attenuation operator ${\cal V}$ might
even be independent of position and, hence, ${\cal L}_{\Delta t}$
has to be computed only once for given medium parameters and 
time step. It can often be done analytically as, for example,
in multiresonant Lorentz models (see Section 10). On the other
hand, by Lemma 7.2, $\|{\cal L}_{\Delta t}\| \leq 1$ for any $\Delta t >0$, 
and one might hope to stabilize the leapfrog scheme by
satisfying the stability condition (\ref{vn}) for ${\cal H}_0$ only, 
that is, $1+{\cal H}_0^2\Delta t^2$ is positive semidefinite. 
This is indeed the case.
The amplification matrix, $\Psi_{t+\Delta t}=
{\cal G}_{\Delta t}\Psi_t$, for the recurrence (\ref{lfnew}),
satisfies the equation
\be
\label{lfam}
{\cal G}_{\Delta t} = {\cal L}_{2\Delta t}{\cal G}_{\Delta t}^{-1}
+ 2\Delta t {\cal L}_{\Delta t}{\cal H}_0\ .
\ee
According to our analysis of the von Neumann stability
condition (\ref{vn}), the anti-Hetmiticity condition
of ${\cal H}_0$ in (\ref{lfnew}) and (\ref{lfam}) 
can be weakened by demanding that ${\cal H}_0$
is related to an anti-Hermitian matrix by a similarity
transformation.
Some important properties of the amplification matrix obtained from (\ref{lfam})
are stated in the following theorem.

{\bf Theorem 9.1}. Suppose there exists a similarity transformation such that
 ${\cal S}^{-1}{\cal H}{\cal S} = {\cal H}_S + {\cal V}_S$ where 
 ${\cal H}_S^*  = -{\cal H}_S$,
the Hermitian part of ${\cal V}_S$ is negative semidefinite, 
${\cal V}_S^{*} + {\cal V}_S \leq 0$, and
${\cal H}_S$ also satisfies the von Neumann stability condition
for the leapfrog scheme, $1 + {\cal H}_S^2\Delta t^2 \geq 0$
(positive semidefinite). Consider 
the amplification matrix ${\cal G}_{\Delta t}$ 
of the modified leapfrog scheme (\ref{lfam}) with
${\cal H}_0 = {\cal S}{\cal H}_S{\cal S}^{-1}$ and
${\cal V}={\cal S}{\cal V}_S{\cal S}^{-1}$. Then
there exists a norm $\|\cdot \|_\mu$ equivalent to
$\|\cdot\|$ such that ${\cal G}_{\Delta t}$
has the following properties:

({\it A}) if $[{\cal H}_0, {\cal V}]=0$,  
\be 
\label{th6a}
\|{\cal G}_{\Delta t}^n\|_\mu \leq 1
\ee
uniformly in $n\geq 0$; 

({\it B}) if $[{\cal H}_0, {\cal V}]\neq 0$, there 
exists a non-negative constant $K_3$ such that 
\be
\label{th6b}
\left\|{\cal G}_{\Delta t} \right\|_\mu
\leq 1 + K_3\Delta t^3\ ,
\ee
for $0<\Delta t<\tau$ and some positive $\tau$.

{\bf Proof}. Part ({\it A}). If ${\cal H}_0$ and ${\cal V}$ commute,
the amplification matrix ${\cal G}_{\Delta t} =
{\cal L}_{\Delta t}{\cal G}_{\Delta t}^0$ 
satisfies (\ref{lfam}), provided ${\cal G}_{\Delta t}^0$
satisfies the same equation for ${\cal V}=0$ (or ${\cal L}_{\Delta t}
=1$), which one can  easily check by substituting  the
solution into (\ref{lfam}). Consider the norm associated with
the similarity transformation ${\cal S}$ of the Hamiltonian,
$\|{\cal A}\|_\mu = \|{\cal S}^{-1}{\cal A}{\cal S}\|$.
The norms $\|\cdot\|_\mu$ and $\|\cdot\|$ are equivalent
(see the proof of Theorem 7.3).  
Since ${\cal H}_0$ satisfies the von Neumann stability condition
and is anti-Hermitian relative to the $\mu$ scalar product, 
$\|{\cal G}_{\Delta t}^0\|_\mu = 
\rho({\cal G}_{\Delta t}^0) = 1$ (according to the analysis
after (\ref{l5})).
By Lemma 7.2, $\|{\cal L}_{\Delta t}\|_\mu\leq 1$,
and we infer that 
$\|{\cal G}_{\Delta t}^n\|_\mu = \|({\cal L}_{\Delta t}{\cal G}_{\Delta
t}^0)^n\|_\mu
\leq \|{\cal L}_{\Delta t}\|^n_\mu \leq 1$
uniformly for $n\geq 0$. 

Part ({\it B}). Solving (\ref{lfam}) by the perturbation theory
in $\Delta t$, it is not hard to find that
\be
{\cal G}_{\Delta t} - {\cal G}_{\Delta t}^V 
= \Delta t^3 {\cal K}_{\Delta t}\ , \ \ \ \ 
{\cal G}_{\Delta t}^V = {\cal L}_{\Delta t/2} 
{\cal G}_{\Delta t}^0{\cal L}_{\Delta t/2}\ ,
\ee
where ${\cal K}_{\Delta t}$ is regular in the vicinity
of $\Delta t=0$ and vanishes whenever ${\cal H}_0$ and
${\cal V}$ commute. On the grid, ${\cal H}_0$
and ${\cal V}$ are bounded operators. Hence we can find
a constant $K_3=\sup_{\Delta t} \|{\cal K}_{\Delta t}\|_\mu$
for some open interval  $0< \Delta t < \tau$. Making use
of the inequality $\|{\cal G}_{\Delta t}^V\|_\mu
\leq \|{\cal L}_{\Delta t/2}\|^2_\mu \leq 1$, we find
\be
\label{th6b1}
\|{\cal G}_{\Delta t}\|_\mu = 
\|{\cal G}_{\Delta t}^V + \Delta t^3 
{\cal K}_{\Delta t}\|_\mu \leq 1 + \Delta t^3 K_3\ ,
\ee
which completes the proof.

The norm deviation of the solution generated
by the modified leapfrog scheme (\ref{lfnew}) from  the stable
solution  generated by ${\cal G}_{\Delta t}^V$ is of order
$O(\Delta t^2)$ for the entire simulation time $T$ and, 
hence, by reducing $\Delta t$ a possible norm
growth can be suppressed as much as desired. Indeed, 
\ba
\nonumber
\left\|{\cal G}_{\Delta t}^n - ({\cal G}_{\Delta t}^V)^n\right\|_\mu&=&
\|\sum_{k=1}^{n-1}{\cal G}_{\Delta t}^{n-k}({\cal G}_{\Delta t}
-{\cal G}_{\Delta t}^V)({\cal G}_{\Delta t}^V)^k\|_\mu\\ \nonumber
&\leq& \Delta t^3 K \sum_{k=1}^{n-1}\|{\cal G}_{\Delta t}^{n-k}\|_\mu\leq
K_3T\Delta t^2 e^{K_3T\Delta t^2} = O(\Delta t^2) \ .
\ea
Since in the continuum limit $\Delta t\rightarrow 0$, 
both the amplification matrices
${\cal G}_{\Delta t}^V$ and ${\cal G}_{\Delta t}$ generate
the same solution and all the powers of the former
are uniformly bounded by construction, a natural question to ask 
is whether one can find a recurrence relation for the function 
$\Psi^V_{n\Delta t}=({\cal G}_{\Delta t}^V)^n\Psi_0$
which could be used in place of (\ref{lfnew}). It is not difficult
to derive an equation for ${\cal G}_{\Delta t}^V$ similar to
(\ref{lfam}), but, unfortunately,
this equation cannot be converted into a simple recurrence relation
for the wave function itself, like 
(\ref{lfnew}), suitable for numerical applications.

It should be noted that if the operator ${\cal L}_{\Delta t}$ in the 
modified leapfrog scheme (\ref{lfnew}) is replaced by another 
${\cal L}_{\Delta t}^s$ such that ${\cal L}_{\Delta t} - 
{\cal L}_{\Delta t}^s = O(\Delta t^3)$ and $\|{\cal L}^s_{\Delta
t}\|_\mu \leq 1$, then the convergence is not violated because
Part B of Theorem 9.1 still holds. This observation
is useful for analytic computation of ${\cal L}_{\Delta t}$.
For example, in the conditions of Theorem 9.1, put ${\cal S}=1$. Let
${\cal V} = {\cal V}_1 + {\cal V}_2$ so that
both ${\cal V}_{1,2}$ have their hermitian parts negative
semidefinite. By using the split (\ref{split}) we get
\be
\label{lvv}
{\cal L}_{\Delta t} = e^{\Delta t {\cal V}} = e^{\Delta t {\cal
V}_1/2}
e^{\Delta t{\cal V}_2}e^{\Delta t{\cal V}_1/2} + O(\Delta t^3)
={\cal L}_{\Delta t}^s + O(\Delta t^3)\ .
\ee
By Lemma 7.2, $\|{\cal L}_{\Delta t}^s\|\leq 1$ for $\Delta t\geq 0$. The operators
${\cal V}_{1,2}$ can be chosen so that their exponentials can be
computed analytically.

\section{Examples of the temporal leapfrog algorithm}
\setcounter{equation}0

There are many possibilities to split the original Hamiltonian ${\cal
H}$ into two parts that satisfy the conditions of Theorem 9.1 and
thereby to make the leapfrog scheme stable and convergent.
Basic guide lines for doing that are as follows. The Hamiltonian
${\cal H}_0$ should contain all the derivative operators in ${\cal H}$
and, yet, the von Neumann condition is easy to establish for ${\cal H}_0$.
It would also be helpful to have an analytic expression for 
${\cal L}_{\Delta t}$ at least up to order $O(\Delta t^3)$.
As an illustration, we discuss multiresonant Lorentz models
and geometric optics. To distinguish between the splits of the 
Hamiltonian in the split and leapfrog algorithms, we shall use
an index $l$ (``leapfrog'') in the latter.

\subsection{Lorentz models}

In the field representation of the Hamiltonian for multiresonant
Lorentz models, we make the following decomposition
\be
\label{lfsplit1}
{\cal H}^F = \pmatrix{{\cal H}_0 & {\cal V}_{FM}\cr {\cal V}_{MF} &0}
+ \pmatrix{0&0\cr 0&{\cal H}_M^F} \equiv {\cal H}_{0l}^F + {\cal V}_l^F
\ee
Thanks to (\ref{lmvmf}) and ${\cal H}_0^* = - {\cal H}_0$, the 
operator ${\cal H}_{0l}^F$ is anti-Hermitian. From (\ref{lmhm}) it 
follows that the Hermitian part of ${\cal V}_l^F$ is negative
semidefinite ($\gamma_a \geq 0$). The exponential of ${\cal V}^F_l$
is easily computed according to (\ref{eHM2}) and (\ref{eHM}).
Let $\xi^a_t$ denote a six-component column whose three upper components
coincide with $\mbf{\xi}^{2a-1}_t$ and three lower components equal 
$\mbf{\xi}_t^{2a}$ (see (\ref{aux})--(\ref{aux2})).
As a result we arrive at the following scheme
\ba
\label{lflm1}
\psi_{t+\Delta t}^F &=& \psi_{t-\Delta t}^F + 2\Delta t {\cal H}_0\psi^F  
+2\Delta t \sum_a {\cal V}_{FMa}\xi^a_t\ ,\\
\label{lflm2}
\xi_{t+\Delta t}^a &=& e^{2\Delta t\, {\cal H}_{Ma}^{F}}\ \xi_{t-\Delta t}^a +
 2\Delta t\ e^{\Delta t\, {\cal H}_{Ma}^{F}}\ {\cal V}_{MFa}\ \psi^{F}_t \ .
\ea
Stability is ensured if ${\cal H}_0^F$ satisfies the von Neumann
condition (\ref{vn}). Eigenvalues of  ${\cal H}_0^F$
satisfy the equation
\be
\label{biswa}
\det(z - {\cal H}_{0l}^F) = z^q \det(z^2 -z{\cal H}_0 - {\cal
 V}_{FM}{\cal V}_{MF}) = 0
\ee
where the non-negative integer $q$ depends on the number of resonances
in the Lorentz model. Non-zero eigenvalues satisfy the so-called
pencil equation whose theory is well developed and might be useful
for more general models \cite{biswa}. Here we shall find a simpler (practical) criterion
sufficient for (\ref{vn}) to hold. Since the plasma frequencies may
 depend
on position, we apply the following general idea \cite{shin}. 
Suppose we have a finite difference scheme with variable 
coefficients in space. Consider a corresponding finite difference scheme
with {\it frozen} coefficients. It is obtained from the original
 scheme by fixing the coefficients to particular values everywhere in space.
A finite difference scheme with variable coefficients is stable
if all the corresponding finite difference
 schemes with frozen coefficients are stable \cite{shin,richt}. So let us fix 
the plasma frequencies to particular values. The spatial dependence
of the eigenfunctions for the pencil problem in (\ref{biswa})
is given by a harmonic factor $\exp(i{\bf k}\cdot{\bf x})$
and the corresponding eigenvalues are $z=\pm i\sqrt{c^2{\bf k}^2
 +\omega_p^2}$,
where $\omega_p^2$ is defined in (\ref{omegap}). Let $k_{max}$
be the maximal norm of all wave vectors of the initial wave packet
and $\omega_p^{max}$ be the maximal value of $\omega_p$ as a function
of position, then a sufficient criterion for stability reads
\be
\label{sclflm}
\Delta t\sqrt{c^2k^2_{max} +(\omega_p^{max})^2} \leq 1\ .
\ee
The scheme (\ref{lfnew}) becomes especially simple in the case of
small attenuation, $\gamma_a < \omega_a$. In the complex
representation of the auxiliary fields (\ref{S}) (cf. (\ref{caux}))
the matter Hamiltonians ${\cal H}_{Ma}^F$ are diagonal and
the action of its exponential is reduced to multiplication by
a complex number $e^{i\nu_a\Delta t}$ (see Section 5).

The stability condition (\ref{sclflm}) can be improved
if one uses the induction representation arranging the 
split according to (\ref{isplit}), that is,
${\cal H}_{0l}^I ={\cal H}_0^I$ and 
${\cal V}^I_l = {\cal V}^I$. In this case the 
conditions of Theorem 9.1 are met if instead of (\ref{sclflm})
we demand a weaker condition
\be
\label{vnh0}
\Delta t ck_{max} \leq 1\ .
\ee
To prove this, we note first that by the similarity transformation
defined in (\ref{main2}) we get
\be
\label{lfsplit2}
{\cal S}^{-1}{\cal H}^I{\cal S}= \pmatrix{{\cal H}_0 & 0\cr 0 & 0} +
\pmatrix{0&{\cal V}_{FM}\cr {\cal V}_{MF}& {\cal H}_M^F} \equiv {\cal H}_S^I
+{\cal V}_S^I\ .
\ee
Then (\ref{vnh0}) is obviously the von Neumann stability condition
for ${\cal H}_S^I$, while the Hermitian part of ${\cal V}_S^I$
is negative semidefinite if $\gamma_a \geq 0$. The scheme
is obtained from (\ref{lfnew}) by replacing $\Psi_t\rightarrow
\Psi^I_t$, ${\cal H}\rightarrow {\cal H}_0^I$ and ${\cal V}
\rightarrow {\cal V}^I$ as defined in (\ref{isplit}). Since
the left hand side of (\ref{lfsplit2}) coincides with ${\cal H}^F$, it can
also be viewed as the leapfrog scheme in the field representation
but with the split different from (\ref{lfsplit1}). This illustrates
the point that the stability condition of the scheme (\ref{lfnew}) depends strongly
on the choice of ${\cal H}_{0l}$.
The price for a
simpler stability condition in the induction representation
is the lack of an explicit form of ${\cal L}_{\Delta t}$.
However, this problem can be circumvented by making use of
(\ref{lvv}).
Indeed, ${\cal V}_S^I = {\cal V}^F + {\cal V}_l^F$ where ${\cal V}^F$
is defined in (\ref{fsplit}). The exponentials
of these operators are computed in Section 5. We also have
$\|\exp(\Delta t {\cal V}^F)\| = 1$ because ${\cal V}^{F*}=
-{\cal V}^F$ and $\|\exp(\Delta t{\cal V}_l^F)\|\leq 1$
by Lemma 7.2 for a non-negative $\Delta t$. 
We set 
\be
\label{llf}
{\cal L}_{\Delta t}^I ={\cal S} e^{\Delta t{\cal V}_l^F/2} 
e^{\Delta t{\cal V}^F}e^{\Delta t{\cal V}_l^F/2}{\cal S}^{-1}
\ee    
so that $\|{\cal L}_{\Delta t}^I\|_\mu =
\|{\cal S}^{-1}{\cal L}_{\Delta t}^I{\cal S}\|\leq 1$.
The operator (\ref{llf}) differs
from ${\cal L}_{\Delta t}$ $= \exp(\Delta t{\cal V}^I)$ $ = 
{\cal S}\exp(\Delta t{\cal V}_S^I){\cal S}^{-1}$ 
by terms of order $O(\Delta t^3)$ and, hence, according to (\ref{lvv}), 
can be used in place of ${\cal L}_{\Delta t}$ in the leapfrog scheme
without destroying its convergence and stability.

\subsection{Geometric optics}

Another simple example is the case of geometric optics. 
For sake of simplicity we assume the medium to have no
magnetic properties. A generalization is straightforward.
Let $\varepsilon = \varepsilon({\bf x})$ be the dielectric
constant of the medium. If the medium is not isotropic,
then $\varepsilon$ is symmetric positive definite $3\times 3$
matrix everywhere in space.
We rewrite Maxwell's equations in the form
\be
\label{schr}
\partial_t \psi_t^I = {\cal H}_G\psi_t^I\ ,\ \ \ \ {\cal H}_G =  
\pmatrix{0 & c\mathbf{\mbf{\nabla}}\times \cr -
c\mathbf{\mbf{\nabla}}\times (\varepsilon^{-1}\  \ )& 0\cr} \ ,
\ee
where the parentheses in $(\varepsilon^{-1})$ 
mean that the induction is first multiplied by $\varepsilon^{-1}$
and then the curl of the resulting vector field is computed. 
Consider the scalar product
\be
\label{sp}
(\psi_1^I,\psi_2^I) = 
\int d{\bf r}\ \psi^{I*}_1\mu \psi_2^I\ ,\ \ \ \mu = 
\pmatrix{\varepsilon^{-1} & 0 \cr 0& 1\cr}\ .
\ee 
In the grid representation of Section 3, the integral is replaced by
the sum over grid points and ${\cal H}_G$ becomes a finite matrix.
The Hamiltonian is anti-Hermitian with respect to this scalar product,
${\cal H}^*_G\mu = -\mu{\cal H}_G$.
Therefore the corresponding $\mu$ norm is preserved in the time evolution
generated by $\exp(t{\cal H}_G)$, that is, 
$(\psi_t^I,\psi_t^I) = 
(\psi_0^I,\psi_0^I)$. The electromagnetic energy of the wave packet
is conserved because it is proportional to the $\mu$ norm of the 
initial state vector. Consequently, we
expect that for a sufficiently small $\Delta t$ the original 
leapfrog scheme (\ref{lf1}), 
\be
\label{lfgo1}
\psi_{t+\Delta t}^I = \psi_{t-\Delta t }^I + 2\Delta t{\cal H}_G\psi_t^{I}\ ,
\ee
becomes stable. To find a sufficient condition for stability, the same
idea of finite difference schemes with frozen coefficients can be used.
It obviously leads to a condition similar to (\ref{vnh0}), 
$$
\Delta t\, c\, k_{max}^\varepsilon\leq  1 \ ,
$$
where $k_{max}^\varepsilon$ is
the maximal norm of all wave vectors in the medium which
can be estimated by $\sqrt{\rho(\varepsilon)}k_{max}$
with $k_{max}$ being the maximal wave vector of the initial
pulse in vacuum. The spectral radius
$\rho(\varepsilon)$ is understood as the maximal spectral
radius of $\varepsilon({\bf x})$ over ${\bf x}$.
If the Fourier basis is used to compute the derivatives,
the algorithm does not violate the Gauss law.
However, the algorithm would not conserve the $\mu$ norm (or energy),
rather a quantity which, in many cases, 
approximates the energy.
Multiplying (\ref{lfgo1}) by $\psi_t^I$ using the scalar product (\ref{sp}),
we infer that
\be
\label{ec}
(\psi_{t+\Delta t}^I, \psi_{t}^I) = (\psi_t^I , \psi_{t-\Delta t}^I)
= \cdots 
= (\psi_{\Delta t}^I, \psi_0^I)\ .   
\ee
By expanding the exponential in $\psi_{t+\Delta t}^I= \exp(\Delta
t{\cal H}_G)\psi_t^I$ into a Taylor series in both sides
of (\ref{ec}) and making use of the 
anti-Hermiticity of ${\cal H}_G$, we find that 
the energy conservation violation is of order $O(\Delta t^2)$.  
Thus, it can be made as small as desired by reducing the time 
step.

\section{Conclusions} 

The initial value problem in Maxwell theory for passive media
has been reformulated in the Hamiltonian formalism. The path integral
representation of the fundamental solution of the Hamiltonian evolution 
equation has been used to develop a time domain numerical algorithm
for solving the initial value problem. The algorithm exhibits the main
advantages of pseudospectral methods for solving differential equations
such as an exponential convergence (and, hence, a greater accuracy), 
the absence of dispersive errors
and numerical efficiency. In addition, the algorithm is unitary,
meaning
that the energy of the initial pulse is conserved whenever the medium
attenuation vanishes (Theorem 6.1). For widely used multiresonant 
Lorentz models,  the algorithm  is unconditionally stable (Theorems 7.1
and 7.3), and, for a generic passive medium, conditional stability
can always be achieved (Theorem 7.4). As the time step $\Delta t$ goes to zero,
the algorithm accuracy is of order $O(\Delta t^2)$ (Theorem
8.1). It is possible to increase the convergence rate (accuracy) up to
any desired order $O(\Delta t^n),\ \ n\geq 2$. However, computational
costs for increasing the accuracy in such a way are not necessarily
lower than those for decreasing the time step in the original
algorithm. An important advantage of the algorithm is that the Gauss
law holds exactly in the process of numerical simulations with no extra
computational cost (Theorem 8.2).

A drawback of the algorithm is related to well known problems 
of the fast Fourier method. Namely, a slower rate of convergence
for non-smooth functions and aliasing. This might, perhaps, limit
the advantages of the algorithm in some type of scattering
problems with complex target geometries. Numerical tests are
needed for a quantitative conclusion.
There are several pseudospectral
methods for approximating the fundamental solution of the Hamiltonian
evolution equation that can help to circumvent this problem. We have 
analyzed one of them and formulated its stability criteria in the case
of general passive media (Theorem
9.1). Numerical tests of the modified leapfrog scheme are presented
in \cite{tests}. The results are compared with known theoretical
and experimental studies of the system investigated (extraordinary
transmission gratings \cite{nature}).
Other methods will be discussed elsewhere as well as the case
when radiation sources (antennas) are included.

It is believed that the proposed algorithm would be useful in numerical
studies of electromagnetic pulse propagation in passive media (e.g., 
foliage, soil, etc), photonic crystals and devices, nonostructured
materials, 
and also in applications to scattering problems
with targets made of dispersive materials.

\vskip 0.2cm
\noindent
{\large\bf Acknowledgments}

I am grateful to Richard Albanese (Brooks, Air Force Base, TX) for suggesting
this project and his continued support. I would like to thank John
Klauder and Tim Olson (University of Florida) for useful
discussions and encouragement. It is my pleasure to express my gratitude to LCAM
(University of Paris-Sud) for the warm hospitality and support of
this project, and my special thanks to Roger Azria and Victor Sidis of
LCAM. I am deeply indebted to Andrei Borisov (LCAM, University of
Paris-Sud) for numerous fruitful discussions and
whose outstanding expertise in computational physics
was invaluable for me.
This work has been supported in part by US Air Force Grants
F4920-03-1-0414 and F49620-01-1-0473.

\section{Appendix}
\subsection{Initial pulse configurations}
\setcounter{equation}0

In principle, any field configuration can serve as the initial
configuration. However 
it is often desired to
have an initial pulse with some specific properties 
(bandwidth, polarization, direction of propagation, etc).
Yet, the initial wave packet should be built of radiation
(propagating)  electromagnetic fields.
A simple method based on the fast Fourier transform algorithm
is given below to obtain initial configurations made
of radiation fields with designated properties.

A general solution of the Maxwell equations in vacuum
can be written in the form
\ba
\label{cc}
{\bf E}_t({\bf r}) &=& \int \, d{\bf k} 
\left( {\bf C}^{(+)}_0( {\bf k})\ e^{i{\bf k}\cdot{\bf r} +ickt} +
 {\bf C}^{(-)}_0({\bf k})\ e^{i{\bf k}\cdot{\bf r} -ickt}\right) 
\equiv { \bf C}^{(+)}_t ({\bf r}) + {\bf C}^{(-)}_t({\bf r})\\
\label{bcc}
{\bf B}_t({\bf r}) &=& \frac{i}{\sqrt{-\Delta}}\ \mbf{\nabla} \times 
\left({ \bf C}^{(+)}_t ({\bf r}) - {\bf C}^{(-)}_t({\bf r})\right)\ .
\ea
The representation (\ref{bcc}) for the magnetic induction follows from
the Maxwell equations and that
the complex amplitudes ${\bf C}_t^{(\pm)}({\bf k})$ 
satisfy the transversality and reality conditions which are,
respectively,
$$
{\bf k}\cdot {\bf C}^{(\pm)}_0 ({\bf k}) = 0\ , \ \ \ \ 
 \overline{\bf C}^{(\pm)}_0 ({\bf k}) ={\bf C}^{(\mp)}_0 (-{\bf k})  \ .
$$
The representation (\ref{cc}) and (\ref{bcc}) 
holds for any moment of time. Hence, we can set $t=0$ to generate suitable initial
conditions for an electromagnetic pulse propagating in empty space by
choosing specific functions ${\bf C}_0^{(\pm)}({\bf r})$.

Consider a few examples. Let there be translational invariance along
the $y$ axis. 
In this case the fields
depend only on $x$ and $z$, i.e., ${\bf r} = (x,0,z)$. Accordingly, the wave vector has the form
${\bf k} = (k_x, 0, k_z)$ and $d{\bf k} = dk_xdk_z$. Let $\hat{\bf e}{}_2 = (0,1,0)$ be the unit vector
along the $y$ axis. Introduce
\be
\label{cc1}
{\bf C}_0^{\pm}({\bf r}) = {\textstyle{\frac 12}}\hat{\bf e}{}_2\ A_0\
e^{-\kappa^2 r^2/2} \ 
e^{\mp i{\bf k}_0\cdot {\bf r}}
 \equiv\hat {\bf e}_2\ C^{\pm}_0 ({\bf r})\ ,
\ee
where $A_0$ and $\kappa$ are real constants, and ${\bf k}_0$ is a fixed wave vector. Note that
the field ${\bf C}_0^{(\pm)}$ is automatically transversal.
The corresponding fields determine suitable initial conditions to generate a pulse propagating 
in the direction of ${\bf k}_0$ whose frequency band is centered at $\omega_0 = ck_0$ and its width
is proportional to $c\kappa$. The pulse
is linearly polarized along the $y$ axis: 
\ba
{\bf E}_0 &=& \hat{\bf e}_2 E_0 \ ,\ \ \ \ \  E_0 = C_0^{(+)} + C_0^{(-)}\ ,\\
{\bf B}_0 &=&  \frac{i}{\sqrt{-\Delta}}\ \mbf{\nabla} \times \hat{\bf
e}{}_2 
\left( C_0^{(+)} - C_0^{(-)}\right)\ .
\ea
The action of the differential operator is defined 
via the fast Fourier transform (see Section 3) on the grid fine enough to support 
the bandwidth limited function (\ref{cc1}). In the Fourier basis,
$i\mbf{\nabla}/\sqrt{-\Delta} \rightarrow -{\bf k}/k$.

To obtain suitable initial conditions for a pulse propagating in the direction ${\bf k}_0$ and 
whose polarization lies in the $xz$-plane, we make use of
the electromagnetic duality of Maxwell theory which states that the dynamics remains
unchanged when electric and magnetic charges switch places and simultaneously 
${\bf E}_t \rightarrow -{\bf B}_t$ and ${\bf B}_t \rightarrow {\bf
E}_t$. 
According to the duality theorem,
we can take
\ba
{\bf E}_0 &=&  \frac{i}{\sqrt{-\Delta}}\ \mbf{\nabla} \times \hat{\bf
e}{}_2 
\left( C_0^{(+)} - C_0^{(-)}\right)\\
{\bf B}_0 & =& -\hat{\bf e}{}_2 \left( C_0^{(+)} + C_0^{(-)}\right)\ .
\ea

Finally, suitable initial conditions for a pulse propagating in the direction ${\bf k}_0$ with a generic
elliptic polarization are obtained by taking a linear combination of the above two initial conditions
for two independent linear polarizations of the pulse. The amplitudes
${\bf C}_0^{(\pm)}({\bf r})$ can also be set numerically from 
actual measurements of a particular pulse of interest.

\subsection{Conductivity and absorbing boundary conditions}

In numerical simulations, the grid in coordinate space 
is of necessity finite. In scattering problems
we are interested in the pulse shape and polarization which we wish to compute in the asymptotically
large coordinate region. This requires that not only the leading edge of the reflected pulse should have
reached the asymptotic region, but also the trailing edge should have
done so as well. 
This is essential if the reflected pulse propagates
in a highly dispersive medium, or the target has a complex shape, or both.  
A complication arises from the very nature of the fast Fourier transform method.
The method is designed to describe periodic functions and, consequently, if the pulse
has a finite amplitude at the edge of the grid, this finite value would appear back at the other edge,
with totally disastrous results for the computation. In quantum computational physics this 
problem is often solved by using an optical potential that absorbs the signal as it reaches the
grid boundary. A similar method can be developed 
for our treatment of Maxwell's theory. 
Before we do so let us point out that an absorbing boundary condition 
is not the only way to solve the problem.
For instance, in the case of a complex target, an ancillary grid may 
be defined in one of the coordinates which extends
to large distances. The pulse may be transferred in a gradual manner 
from the small grid (near the target)
to this larger grid to prevent the pulse from  ever reaching the edge of the small grid. This technique
can also be applied to generate a pulse by an antenna of a complex construction. The dynamics of the 
portion of the pulse on the larger grid may be treated analytically 
(if dispersion properties of the medium 
are not too complex).

 In quantum mechanics absorbing boundary conditions are made by adding an imaginary potential
to the Hamiltonian with support near the grid edges. In the Maxwell
theory, 
the same can be achieved by adding
conductivity which gradually increases as the grid edges are approached. An interaction of conducting
media with electromagnetic radiation is described by  Ohm's law,
\be
{\bf J}_t =\sigma {\bf E}_t\ ,
\ee
combined with Maxwell's equation (\ref{m1}), where the displacement
current is amended as $\partial_t{\bf D}_t\rightarrow
\partial_t{\bf D}_t + (4\pi \sigma/c){\bf E}_t$ with 
$\sigma = \sigma({\bf r}) $ being the conductivity of the medium.
Consider a linearly polarized plane wave moving along the $z$ axis. Let $\tilde{ E}_\omega$ be 
the Fourier transform of the only component of the electric field
${E}_t$. 
Disregarding for a moment any possible anomalous dispersion of 
the medium, we find that $\tilde{ E}_\omega$ satisfies the equation
\be
\label{con2} 
 \partial_z^2\tilde{E}_\omega(z) +\left[\frac{\omega^2}{c^2} 
- \frac{4\pi i \omega}{c}\, \sigma(z)\right]\, \tilde{E}_\omega(z)
=0\ .
\ee  
Equation (\ref{con2}) is identical to the stationary Schroedinger 
equation with an optical (absorbing) potential 
being proportional to
$\sigma(z)$.  In simulations of quantum wave packets it has been 
found that one of  the optimal potentials has 
the form \cite{cond}
\be
\label{con3}
\sigma(z) = (n+1)^{-1}\, \sigma_n\ (z/L)^n\ , \ \ \ n\geq 2\ ,
\ee
in the interval $z\in [0,L]$ and $\sigma (z)=0$ otherwise. So, our next task is to find an optimal
 constant $\sigma_n$ such that that the conducting layer would not reflect or transmit electromagnetic
energy in some designated frequency band.
Maxwell's equations in a conducting medium are form-invariant under
the scaling transformations
\be
\omega \rightarrow \beta\omega\ ,\ \ \ 
\sigma({\bf r}) \rightarrow \beta\sigma(\beta {\bf r})\ ,\ \ \ \ 
\tilde{\psi}^F_\omega({\bf r}) \rightarrow \alpha\tilde{\psi}^F_\omega(\beta {\bf r})\ ,
\ee
where $\alpha$ and $\beta$ are positive constants.
If the conductivity $\sigma({\bf r})$ was found optimal for a
frequency 
$\omega$ and over a length $L$, then 
the optimal conductivity for a frequency $\beta\omega$ would be $\beta\sigma(\beta {\bf r})$ and the 
new length over which it is taken to act would be $L/\beta$.

Let us first study
the reflectivity of the absorbing layer. Suppose, $\sigma(z) = \sigma_0\theta(z)$ where
$\theta(z)$ is the Heaviside function. If a monochromatic linearly polarized wave coming from
the negative $z$ region has an amplitude one, then the reflected wave
has the amplitude \cite{landau}
\be
\label{abs1}
\eta_\omega = \frac{1-\nu_\omega}{1+\nu_\omega}\ ,\ \ \ \ \
\nu_\omega^2 = 
1 +\frac{4\pi i \sigma_0}{\omega} 
\equiv 1 + i q_\omega\ .
\ee
The energy of the reflected wave is 
\be
\label{abs2}
R_\omega = |\eta_\omega |^2 = q^2_\omega/4 + O(q^4_\omega)\ .
\ee
$R_\omega$ increases as the ratio $q_\omega$ gets higher and  is small if $q^2_\omega /4 <\!\!< 1$.
Consider now $\sigma(z)$ which monotonically increases from $z=0$ 
in the positive $z$ direction. Let $k_\omega
=k_\omega(z)$
be a local wave vector of the wave in the conducting medium,
\be
\label{abs3}
k_\omega (z) = c^{-1}\omega\sqrt{1 + i q_\omega(z)}\ .
\ee
We shall argue that the reflection is negligibly 
small if the local wave vector does not significantly changes
over a distance of order $k_\omega^{-1}$, that is,
\be
\label{abs4} 
\left\vert\frac{k_\omega(z+\delta z) - k_\omega(z)}{k_\omega(z)}\right\vert <\!\!< 1\ ,\ \ \ \ \ 
\delta z \sim k_\omega^{-1}\ .
\ee
Making a linear approximation in (\ref{abs4}), we infer that
\be
\label{abs5}
|\partial_z q_\omega(z)| <\!\!< (2c)^{-1}\,\omega\, \left(1 + q_\omega^2(z)\right)^{3/4}\ ,
\ee
which must be valid for all values of $q_\omega(z)$ including small ones when the reflection is small.
Inequality (\ref{abs5}) allows us to reverse the argument, that is,
the reflection is small if $|\partial_z q_\omega(z)| <\!\!< (2c)^{-1}\,\omega$. Let $\sigma_L$ be an
average conductivity over a layer of width $L$, $\sigma_L = L^{-1}\int_0^L dz \sigma(z)$. In particular,
for (\ref{con3}), $\sigma_L = \sigma_n$. For a monotonically increasing function, the derivative
can be approximated as $|\partial_z q_\omega(z)| \approx \sigma_L/L$. This leads to a necessary 
condition on conductivity to suppress the reflection, namely,
\be
\label{abs6}
\sigma_L < \frac{\omega^2 L}{8\pi c}\ .
\ee
Our analysis is valid if the higher derivatives of $\sigma(z)$ are not
large. 
This condition requires that
the exponent
$n$ in (\ref{con3}) should not be less than two to insure a smooth behavior at $z=0$.

The transmission can be estimated as follows. Suppose the pulse occupies a compact region $\Omega$.
Let  ${\cal E}_t^\Omega$ be the pulse energy.
The pulse looses its energy as it propagates through a conducting medium according to Ohm's law, so
that 
\be
\label{abs7}
c^{-1}\partial_t  {\cal E}_t^\Omega = - (2c)^{-1} \int_\Omega d{\bf
  r}\ \sigma\ {\bf E}_t^2 \leq -8\pi \sigma_\Omega {\cal E}_t^\Omega\ ,
\ee
where $ \sigma_\Omega = \max_\Omega \sigma$.
Therefore the pulse energy decay can be bounded from above by
\be
\label{abs9}
{\cal E}_t^\Omega \leq e^{-8\pi t\sigma_\Omega} \, {\cal E}^\Omega_0\ .
\ee
In the one dimensional case (\ref{con3}), $\sigma_\Omega =\sigma_L/(n+1)$. For the time
$t=L/c$ needed for a pulse to get through the layer of width $L$, 
the attenuation should be large, that is,
$8\pi L\sigma_L/c(n+1) >\!\!> 1$. Thus, the necessary conditions to suppress both  transmission and
reflection (that is, to ensure an almost total absorption) of the pulse are
\be
\label{abs10}
\frac{(n+1)c}{8\pi L} < \sigma_L < \frac{\omega^2 L}{8\pi c}\ .
\ee

By changing the Hamiltonian ${\cal H}^Q$, the conducting layer can be
included
into the split or leapfrog algorithm. Since the conducting layer 
produces attenuation, the conductivity $\sigma$ must be included
into the operator ${\cal L}_{\Delta t}$ in the modified leapfrog
scheme. It is also possible to create an absorbing and non-reflecting
layer by using a passive medium (e.g. a Lorentz model). The analysis
of the medium properties would be similar to that for a conducting
layer. In fact, using a layer of a passive medium would offer more
flexibility in solving the grid boundary problem.

\end{document}